\documentstyle[amsfonts,leqno]{article}
\let\fraktur\frak

\def\colon{{:}\;}

\newcommand{\sgroup}[2]{\raisebox{-2ex}
{\begin{picture}(26,24)(0,0)
\put(3,15){\vector(1,0){20}}
\put(3,10){\vector(1,0){20}}
\put(10,3){${\scriptstyle #1}$}
\put(10,18){${\scriptstyle #2}$}
\end{picture}}}

\newcounter{th}

\newtheorem{th}[th]{Th\'{e}or\`{e}me}
\newtheorem{prop}[th]{Proposition}
\newtheorem{cor}[th]{Corollaire}
\newtheorem{lem}[th]{Lemme}
\newtheorem{exem}[th]{Exemple}
\newtheorem{exems}[th]{Exemples}

\title{INT\'{E}GRATION SYMPLECTIQUE
DES VARI\'{E}T\'{E}S DE POISSON R\'{E}GULI\`{E}RES}
\author{F. ALCALDE CUESTA\thanks{Recherche support\'ee par
D.G.I.C.Y.T. Espagne (Proyecto PB90-0765) et Xunta de Galicia
(Proxecto XUGA20704B90)}
et G. HECTOR}
\date{}
\begin{document}
\maketitle

\begin{abstract}
A {\bf symplectic integration} of a Poisson manifold $(M,\Lambda)$
is a symplectic groupoid
$(\Gamma,\eta)$ which {\bf realizes}
the given Poisson  manifold, i.e. such that the space of units
$\Gamma_0$ with the induced Poisson  structure $\Lambda_0$ is isomorphic to
$(M,\Lambda)$.  This notion was introduced by A. Weinstein in
\cite{W2} in order to quantize Poisson manifolds by quantizing their symplectic
integration.
Any Poisson manifold can be integrated by a {\bf local} symplectic groupoid
(\cite{CDW}, \cite{Ka}) but already for regular Poisson manifolds there are
obstructions to global integrability (\cite{AM}, \cite{D2}, \cite{He},
\cite{M}, \cite{W2}).

The aim of this paper is to summarize all the known obstructions and present a
sufficient
topological condition for integrability of regular Poisson manifolds; we will
indeed describe a concrete procedure for this integration. Further our
criterion will provide necessary and sufficient if we require $\Gamma$ to be
Hausdorff, which is a suitable condition to proceed to Weinstein's program of
quantization.
These integrability results may be interpreted as an generalization
of the Cartan-Smith proof of Lie's third theorem in the infinite dimensional
case.
\end{abstract}

\section{Introduction}

\hspace{1em}
Un {\bf groupo\"{\i}de symplectique} $(\Gamma,\eta)$ est un
groupo\"{\i}de de Lie muni d'une structure symplectique $\eta$
compatible avec la multiplication partielle
(\cite{CDW}, \cite{Ka}, \cite{W2}). L'espace
des unit\'{e}s $\Gamma_{0}$ est muni d'une structure de Poisson
canonique $\Lambda_{0}$ pour laquelle la projection source $\alpha$ est
un morphisme de Poisson. Le probl\`{e}me
de l'{\bf int\'{e}gration symplectique} d'une vari\'{e}t\'{e} de
Poisson $(M,\Lambda)$ a \'{e}t\'{e} pos\'{e} par A. Weinstein dans \cite{W2};
il consiste \`{a} construire un groupo\"{\i}de symplectique
$(\Gamma,\eta)$ dont l'espace des unit\'es $(\Gamma_0,\Lambda_0)$ est isomorphe
\`a  $(M,\Lambda)$. Quitte \`{a} d\'{e}ployer $\Gamma$ (voir \cite{Pr1}), on
peut supposer que les fibres de $\alpha$ sont connexes et simplement connexes;
une telle int\'{e}gration symplectique est dite {\bf universelle}.
L'id\'{e}e de A. Weinstein est d'utiliser l'int\'{e}gration symplectique pour
quantifier les vari\'{e}t\'{e}s de Poisson
(voir \cite{W2}, \cite{W4}, \cite{WX}).
\vspace{2ex}

L'int\'{e}gration par un groupo\"{\i}de symplectique {\bf local}
a \'{e}t\'{e} r\'{e}alis\'{e}e
dans \cite{CDW} et \cite{Ka}. Pour l'int\'{e}gration globale,
il est naturel de s'int\'{e}resser d'abord aux vari\'{e}t\'{e}s de
Poisson {\bf r\'{e}guli\`{e}res} (i.e. celles dont le
feuilletage caract\'{e}ristique $\cal F$ est r\'{e}gulier) et plus
pr\'{e}cisement aux vari\'{e}t\'{e}s de Poisson r\'{e}guli\`{e}res
dont tout cycle \'{e}vanouissant est trivial. Cette condition situe le
probl\`{e}me dans la cat\'{e}gorie des vari\'{e}t\'{e}s de Poisson
s\'{e}par\'{e}es ce qui permet des calculs cohomologiques.
Dans ce contexte, on a les r\'{e}sultats suivants:

\noindent
1) il existe des vari\'{e}t\'{e}s de Poisson non int\'{e}grables:
P. Dazord construit une obstruction \`{a}
l'int\'{e}grabilit\'{e} de certaines vari\'et\'es de Poisson
dans \cite{D1} et \cite{D2};

\noindent
2) les vari\'{e}t\'{e}s de Poisson {\bf totalement asph\'{e}riques} (i.e.
celles pour lesquelles le $\pi_{2}$ des feuilles de $\cal F$ est nul)
sont int\'{e}grables (\cite{DH}).
\vspace{2ex}

Le premier but de ce travail est d'exhiber une condition
topologique suffi\-sante pour int\'{e}grer les
vari\'{e}t\'{e}s de Poisson r\'{e}guli\`{e}res sans cycle \'{e}vanouissant qui
g\'en\'eralise et r\'esume tous les r\'esultats d'int\'egration connus
\`a ce jour. En fait, on donne une construction explicite
de leur int\'egration symplectique universelle.

Les vari\'et\'es de Poisson r\'eguli\`eres dont l'int\'egration symplectique
est s\'epar\'ee forment une bonne cat\'egorie pour
la pr\'equantification propos\'ee par A. Weinstein.
Pour ces vari\'et\'es, on obtient une condition n\'ecessaire et suffisante
d'int\'egration.

La construction de l'int\'egration s'appuie sur une th\'eorie de
{\bf fibr\'es principaux \`a groupo\"{\i}de  structural}
qui remplace l'approche par les
{\bf r\'ealisations isotropes de Libermann} d\^ue \`a  P. Dazord  (\cite{D1},
\cite{D2}).
On est donc amen\'e \`a introduire dans ce contexte les notions
famili\`eres pour les fibr\'es principaux classiques:
{\bf cocycle, connexion, courbure et classe de Chern}.
Cela permet d'expliciter l'analogie avec les r\'esultats d'int\'egrabilit\'e de
B. Kostant (\cite{Ks}), J. M. Souriau (\cite{So}) et A. Weil (\cite{W}).
Cette approche am\`ene naturellement \`a un autre aspect:
l'utilisation de techniques cohomologiques pour
"int\'egrer" les invariants cohomologiques de la structure de Poisson.

La partie technique du travail consiste \`a d\'emontrer un th\'eor\`eme
d'annulation
pour la cohomologie (\`a valeurs dans un faisceau)
d'une submersion; celui-ci est le coeur de la preuve du
th\'eor\`eme d'int\'egration.
\vspace{1ex}

\subsection{R\'{e}sultats}

\hspace{1em}
D'apr\`es \cite{DH}, une structure de Poisson r\'{e}guli\`{e}re
$\Lambda$ sur une vari\'{e}t\'{e} $M$ est d\'{e}termin\'{e}e par
le {\bf feuilletage caract\'{e}ristique} $\cal F$ engendr\'{e}
par les champs hamiltoniens $X_{f}$ et la {\bf forme feuillet\'{e}e
symplectique} d\'efinie par
$ \sigma(X_{f},X_{g}) = \{f,g\} $.
C'est ce point de vue
{\bf feuillet\'e} qui va permettre d'introduire les objets essentiels dans la
construction de l'int\'egration symplectique.
\vspace{2ex}

Toute sph\`{e}re tangente \`{a} $\cal F$ peut \^{e}tre d\'{e}form\'{e}e
transversalement en une famille continue $D$ de sph\`{e}res tangentes (voir
\S\ref{stdif}). L'int\'{e}gration de $\sigma$ sur ces sph\`{e}res
d\'{e}finit une {\bf fonction d'aire} sur une transversale \`{a} $\cal F$.
Si cette fonction  poss\`{e}de un point critique, on dira que la {\bf
d\'{e}formation transverse} $D$ est {\bf symplectiquement \'{e}vanouissante};
une telle d\'{e}formation sera dite {\bf triviale} si la fonction d'aire est
constante. Les diff\'{e}rentielles des fonctions d'aire engendrent
un sous-groupo\"{\i}de ${\cal P}(\Lambda)$ du fibr\'e conormal
\( \nu^{\ast}{\cal F} \)
appel\'{e} le {\bf groupo\"{\i}de des p\'{e}riodes
sph\'{e}riques d\'{e}riv\'{e}es} (ou simplement le
{\bf groupo\"{\i}de d\'eriv\'e}) de $\Lambda$.

Comme dans les probl\`{e}mes de l'int\'{e}gration des groupes locaux
(\cite{V}) et des suites d'Atiyah (\cite{AM}, \cite{M}) ou dans le
probl\`eme de la  pr\'{e}quantification des vari\'{e}t\'{e}s symplectiques
(\cite{Ks}, \cite{So}, \cite{W}),
l'int\'{e}grabilit\'{e} des vari\'{e}t\'{e}s de Poisson sera
carac\-t\'{e}ris\'{e}e par leurs p\'{e}riodes (sph\'eriques d\'eriv\'ees).
Celles-ci doivent \^etre annul\'ees
si l'on veut faire dispara\^{\i}tre les obstructions cohomologiques \`a
l'int\'egration. Par cons\'equent, le groupo\"{\i}de quotient $\cal G$ du
fibr\'e conormal
\( \nu^{\ast}{\cal F} \)
par le groupo\"{\i}de d\'eriv\'e
\( {\cal P} (\Lambda) \)
sera l'ingr\'edient fondamental de l'int\'egration symplectique.
On r\'esumera cette id\'ee
en dissant que $\cal G$ est le {\bf groupo\"{\i}de structural} de
$\Lambda$.
\vspace{1ex}

\begin{th}         \label{th1}
Soit $(M,\Lambda)$ une vari\'et\'e de Poisson r\'eguli\`ere.
Le groupo\"{\i}de structural $\cal G$ est un groupo\"{\i}de de Lie
si et seulement si

\noindent
i) le groupo\"{\i}de d\'{e}riv\'{e} ${\cal P}(\Lambda)$ est un
sous-groupo\"{\i}de de Lie de
\( \nu^{\ast}{\cal F} \)
\'{e}tal\'{e} sur $M$,
i.e. toute d\'{e}formation symplectiquement \'{e}vanouissante est
triviale;

\noindent
ii) le groupo\"{\i}de d\'{e}riv\'{e} ${\cal P}(\Lambda)$ est plong\'{e} dans
\( \nu^{\ast}{\cal F} \).

\noindent
En outre, $\cal G$ est s\'{e}par\'{e} si et seulement si
${\cal P}(\Lambda)$ est ferm\'{e}.
\end{th}
\vspace{1ex}

La condition (i) est l'obstruction \`a l'integrabilit\'e mise en
lumi\`ere par A. Weins\-tein dans
\cite{W2}. Les conditions (i) et (ii) impliquent que le groupo\"{\i}de
d\'{e}riv\'{e} est un {\bf r\'{e}seau} de $(M,\Lambda)$ au sens
de \cite{D1} et \cite{D2}; il s'agit de l'obstruction de P. Dazord.
En outre, ces deux conditions entra\^{\i}nent que les fibres de
${\cal P}(\Lambda)$ sont
ferm\'{e}es discr\`{e}tes; en restriction aux feuilles de $\cal F$,
on retrouve ainsi
l'obstruction \`a l'int\'egrabilit\'e des alg\'ebro\"{\i}des de Lie transitifs
de R. Almeida et P. Molino (\cite{AM}) et K. Mackenzie (\cite{M}).
Bref, toutes ces
obstructions \`a l'int\'egrabilit\'e n'interviennent en fait que pour assurer
que le groupo\"{\i}de structural $\cal G$ est un groupo\"{\i}de de Lie.
\vspace{1ex}

\begin{exems}   \label{exobs}
{\em
(1) Soit $\cal F$ le feuilletage horizontal sur
\( M = S^{2} \times {\Bbb R} \)
d\'efini par l'\'equation $dt=0$.
Soit $v_{0}$ la forme volume canonique sur
$S^{2}$. La forme feuillet\'{e}e
symplectique $\sigma$ repr\'{e}sent\'{e}e par
\( \omega = (1+t^{2}) v_{0} \)
d\'efinit une structure de Poisson $\Lambda = ({\cal F},\sigma)$ sur $M$.
Le groupo\"{\i}de d\'{e}riv\'{e}, qui est engendr\'e par la 1-forme
$2tdt$, n'est pas un groupo\"{\i}de de Lie.
\vspace{1ex}

Cet exemple de A. Weinstein (voir \cite{W2}) montre que le groupo\"{\i}de
structural $\cal G$ n'est pas en g\'{e}n\'{e}ral un groupo\"{\i}de de Lie.
L'exemple suivant montre que
la condition de plongement est aussi essentielle:
\vspace{1ex}

\noindent
(2) La structure de Poisson sur $S^2 \times {\Bbb R}$ de
l'exemple (1) induit une
structure de Poisson $\Lambda$ sur l'ouvert $M$
obtenu en \^otant le point $(N,0)$, o\`u $N$ est le p\^ole nord de $S^2$.
Il n'y a plus de d\'eformation symplectiquement \'evanouissante,
mais le groupo\"{\i}de structural
$\cal G$ n'est toujours pas un groupo\"{\i}de de Lie.
\vspace{1ex}

Si le groupo\"{\i}de d\'eriv\'e est plong\'{e} dans
\( \nu^{\ast}{\cal F} \)
\'{e}tal\'{e} sur $M$, alors les fibres sont discr\`{e}tes,
mais ces deux conditions ne sont pas \'{e}quivalentes d'apr\`es
l'exemple (2).}
\end{exems}
\vspace{1ex}

Les analogies \'etablies et le proc\'ed\'e de construction de l'int\'egration
symplectique (voir \S4) permettent d'interpr\'eter le th\'eor\`eme suivant
comme une extension du th\'{e}or\`{e}me de Cartan-Smith qui
\'etablit le 3\`{e}me th\'{e}or\`{e}me de Lie (cf. \cite{V}):
\vspace{1ex}

\begin{th}    \label{th2}
Soit $(M,\Lambda)$ une vari\'{e}t\'{e} de Poisson r\'{e}guli\`{e}re dont
tout cycle \'{e}vanouissant est trivial. Si

\noindent
i) le groupo\"{\i}de structural $\cal G$ est un groupo\"{\i}de de Lie,

\noindent
ii) la vari\'et\'e de Poisson $(M,\Lambda)$ est
{\bf transversalement compl\`ete}, i.e. il existe un suppl\`ementaire
du fibr\'e tangent du feuilletage image r\'eciproque de $\cal F$ sur
$\cal G$,

\noindent
alors $(M,\Lambda)$ est int\'egrable.
\end{th}
\vspace{1ex}

Si $\cal G$ est s\'epar\'e, la vari\'et\'e de
Poisson $(M,\Lambda)$ est \'evidemment transversalement compl\`ete
(voir \S\ref{comp}).
D'autre part, le groupo\"{\i}de structural d'une vari\'et\'e de Poisson
{\bf totalement asph\'erique} (i.e. tout cycle \'evanouissant
est trivial et le $\pi_2$ des feuilles est nul) est isomorphe \`a
$\nu^\ast{\cal F}$ et l'on retrouve le th\'eor\`eme
fondamental de \cite{DH}:
\vspace{1ex}

\begin{cor}   \label{totasp}
Si $(M,\Lambda)$ est une vari\'et\'e de Poisson r\'eguli\`ere
totalement asph\'e\-rique, alors $(M,\Lambda)$ est int\'egrable. $\Box$
\end{cor}
\vspace{1ex}

R\'eciproquement, on montre le th\'eor\`eme d'obstruction suivant qui r\'esume
les obstructions pr\'ealablement \'etablies (\cite{AM}, \cite{D2}, \cite{M},
\cite{W2}) dans le
cas g\'en\'eral des vari\'et\'es de Poisson r\'eguli\`eres:
\vspace{1ex}

\begin{th}   \label{th3}
Si une vari\'{e}t\'{e} de Poisson r\'eguli\`ere
$(M,\Lambda)$ est int\'{e}grable,
alors le groupo\"{\i}de structural $\cal G$ est un groupo\"{\i}de de Lie.
\end{th}
\vspace{1ex}

D'apr\`es le th\'eor\`eme~\ref{th2}, il semble raisonnable d'ajouter une
nouvelle obstruction, \`a savoir la compl\'etion
transverse de $(M,\Lambda)$. En fait, le proc\'ed\'e de \cite{He} permet de
montrer une condition n\'ecessaire proche de la compl\'etion transverse (en
\'elargissant la notion de {\bf cycle \'evanouissant coh\'erent} de
\cite{He}). N\'eanmoins, le probl\`eme de l'existence ou non d'un
suppl\'ementaire
invariant reste ouvert.
\vspace{2ex}

La quantification d'apr\`{e}s A. Weinstein (voir \cite{W2} et \cite{WX})
am\`{e}ne \`{a} s'int\'{e}resser aux vari\'{e}t\'{e}s de Poisson dont
l'int\'{e}gration symplectique est s\'{e}par\'{e}e:
\vspace{1ex}

\begin{th}   \label{th4}
Une vari\'{e}t\'{e} de Poisson r\'{e}guli\`{e}re $(M,\Lambda)$ est
int\'{e}grable par un
groupo\"{\i}de symplectique s\'{e}par\'{e} si et seulement si tout cycle
\'evanouissant est
trivial et le groupo\"{\i}de structural $\cal G$ est un groupo\"{\i}de de
Lie s\'epar\'e.
$\Box$
\end{th}
\vspace{1ex}

\hspace{-4pt}
Les conditions d'int\'egration se simplifient si le {\bf groupo\"{\i}de
d'homotopie} $\Pi_{1}({\cal F})$
(voir \S\ref{intP}) est localement trivial.
On remarque que les cycles \'evanouissants d'un tel feuilletage sont
triviaux. La proposition~\ref{pi2} permet d'\'enoncer le r\'esultat suivant:
\vspace{1ex}

\begin{cor}   \label{cor}
Une vari\'{e}t\'{e} de Poisson r\'{e}guli\`{e}re $(M,\!\Lambda)$ \`{a}
groupo\"{\i}de d'homotopie localement trivial est int\'{e}grable si et
seulement si le groupo\"{\i}de d\'{e}riv\'{e} ${\cal P}(\Lambda)$
est un sous-groupo\"{\i}de de Lie de
\( \nu^{\ast}{\cal F} \)
\`{a} fibres ferm\'{e}es discr\`{e}tes. $\Box$
\end{cor}

\begin{exems}
{\em On donne des exemples de structures de Poisson $\Lambda=({\cal F},\sigma)$
dont le groupo\"{\i}de d'homotopie est localement trivial:

\noindent
1) Si $\cal F$ est d\'{e}fini par une action localement
libre d'un groupe de Lie $G$ et puisque $G$ est asph\'{e}rique,
$\Lambda$ est int\'egrable d'apr\`es le corollaire~\ref{totasp}.

\noindent
2) Si $\cal F$ est {\bf riemannien complet}, $\Lambda$ v\'erifie la
condition du corollaire~\ref{cor} d'apr\`es \cite{A1}.

\noindent
3) Une {\bf structure cosymplectique} (\cite{L}) sur une vari\'et\'e $M$
est la donn\'ee d'une
1-forme ferm\'ee $\theta$ et d'une 2-forme ferm\'ee $\omega$ telles que
$\theta \wedge \omega^k$ soit un forme volume sur $M$. Le feuilletage $\cal F$
d\'efini par l'\'equation $\theta=0$ et la forme
feuillet\'ee $\sigma$ repr\'esent\'ee par $\omega$ d\'efinissent une
structure de Poisson $\Lambda$ sur $M$. Si le champ de Reeb $R$ (d\'efini par
$i_R \theta = 1$ et $i_R \omega =0$) est complet, le $\Bbb R$-feuilletage de
Lie
$\cal F$
est complet et donc v\'erifie la propri\'et\'e du corollaire~\ref{cor}. Puisque
le groupo\"{\i}de d\'eriv\'e est nul, $\Lambda$ est int\'egrable. En fait,
puisque
$\omega$ est ferm\'ee, $\Lambda$ est toujours int\'egrable d'apr\`es \cite{DH}.

\noindent
4) Soit $\cal F$ un feuilletage {\bf totalement g\'{e}od\'{e}sible}. En
proc\'edant comme dans \cite{A1}, ce cas se ram\`ene au cas (1) \`a l'aide du
th\'eor\`eme de structure de G. Cairns (voir \cite{Mo}).}
\end{exems}

\begin{exems}
{\em Soient $\cal F$ un feuilletage transversalement orientable
de codimension 1 et $\Lambda=({\cal F},\sigma)$ une structure de Poisson
sur une vari\'et\'e compacte $M$ de dimension 3.
D'apr\`es le th\'eor\`eme de stabilit\'e de Reeb et le th\'eor\`eme de Novikov,
il y a trois cas possibles:

\noindent
1) $\cal F$ est totalement asph\'erique et donc $\Lambda$ est int\'egrable
d'apr\`es le corollaire~\ref{totasp};

\noindent
2) $\cal F$ est la fibration triviale en sph\`eres et
$\Lambda$ est int\'egrable si la fonction d'aire est
constante d'apr\`es le corollaire~\ref{cor};

\noindent
3) $\cal F$ poss\`ede une composante de Reeb qui supporte un cycle
\'evanouissant non trivial et dans ce cas aucune structure de Poisson $\Lambda$
n'est int\'egrable d'apr\`es \cite{He}.}
\end{exems}

\setcounter{th}{0}
\setcounter{equation}{0}

\section{Vari\'et\'es de Poisson r\'eguli\`eres}

\hspace{1em}
Le but de cette section est de rappeler la description des vari\'{e}t\'{e}s
de Poisson r\'{e}guli\`{e}res en termes de formes feuillet\'{e}es
(voir \cite{DH} et \cite{He}).
\vspace{1ex}

\subsection{Formes feuillet\'{e}es et formes pures}

i) Soit $(M,{\cal F})$ une vari\'{e}t\'{e} feuillet\'{e}e.
Soient
\( (\Omega^{\ast}(M),d) \)
le complexe de De Rham et
\( (\Omega^{\ast}M,{\cal F}),d) \)
le sous-complexe des
{\bf formes relatives}, i.e. des formes qui
s'annulent sur les feuilles de $\cal F$.
Le complexe
\( (\Omega^{\ast}({\cal F}),d_{\cal F}) \)
des
{\bf formes feuillet\'{e}es}
est le complexe quotient et sa cohomologie
\( H^{\ast}({\cal F}) \)
est la
{\bf cohomologie feuillet\'{e}e de}
$(M,{\cal F})$.
\vspace{2ex}

\noindent
ii) Le choix d'un suppl\'{e}mentaire $N{\cal F}$ de $T{\cal F}$
permet d'obtenir de bons repr\'{e}sen\-tants
des formes feuillet\'{e}es. Soit
\( \nu^{\ast}{\cal F} \)
le fibr\'{e} conormal dont les sections sont les 1-formes relatives.
La d\'{e}composition
\( T^{\ast}M = \nu^{\ast}{\cal F} \oplus T^{\ast}{\cal F} \)
induit une d\'{e}composition de $\Omega^{\ast}(M)$ et $d$ en
{\bf formes pures} et {\bf composantes pures}
\begin{eqnarray*}
\Omega^{r}(M) = \raisebox{-1ex}{$\stackrel{\bigoplus}{\scriptstyle p+q=r}$}
\Omega^{p,q}(M)
\;\;\;\;\;\;\;
d = d_{0,1} + d_{1,0} + d_{2,-1}
\vspace{1ex}
\end{eqnarray*}

\noindent
iii) Si l'on consid\`ere
\( \Omega^{r}(M) = \bigoplus \Omega^{p,q}(M) \)
comme un module filtr\'e de degr\'e filtrant $p$, on obtient
la {\bf suite spectrale de Leray-Serre}
$$
E_{2}^{p,q}({\cal F}) \Longrightarrow H^{p+q}(M)
$$
Puisque le terme
\( E_{0}^{p,q}({\cal F}) = \Omega^{p,q}(M) \)
et la diff\'{e}rentielle $d_{0}$ est induite par $d_{0,1}$, on a:
$$
E_{1}^{0,q}({\cal F}) \cong H^{q}({\cal F})
$$
iv) Soit $L$ est la forme de Liouville sur
\( \nu^{\ast}{\cal F} \).
Les champs $H$ d\'efinis par

\begin{equation} \label{rel}
i_{\textstyle H} L = i_{\textstyle H} dL = 0
\end{equation}

\noindent
engendrent un feuilletage $\cal H$
appel\'e le {\bf feuilletage relev\'e de $\cal F$}.
Le fibr\'e tangent $T{\cal H}$ est le sous-fibr\'e horizontal
de la {\bf connexion partielle de Bott}
d\'{e}finie par
$$
\nabla_X \mu = i_X d\mu
$$
pour tout champ
\( X \in {\fraktur X}({\cal F}) \)
et toute section $\mu$ de
\( \nu^{\ast}{\cal F} \).
Bref,
\( \nu^{\ast}{\cal F} \)
est un fibr\'e vectoriel {\bf $\cal F$-feuillet\'e} (au sens de \cite{Mo})
dont les {\bf sections feuillet\'ees}
sont les {\bf 1-formes basiques}, i.e. les 1-formes $\mu$ telles que
\( i_X \mu = i_X d\mu = 0 \)
pour tout
$X \in {\fraktur X}({\cal F})$.
\vspace{2ex}

\noindent
v) L'espace
\( \Omega^{q}({\cal F};\nu^{\ast}{\cal F}) \)
des
{\bf q-formes feuillet\'{e}es \`{a} valeurs dans
\( \nu^{\ast}{\cal F} \)}
est l'espace des sections du fibr\'{e} vectoriel
\( (\bigwedge^{q} T^{\ast}{\cal F}) \otimes \nu^{\ast}{\cal F} \).
La diff\'{e}rentielle ext\'{e}rieure covariante
$$
\begin{array}{ll}
d_{\cal F} \omega(X_{1}, \dots ,X_{q+1}) = &
\sum_{i=1}^{q+1} (-1)^{i+1} \nabla_{\textstyle X_{i}}
\omega(X_{1}, \dots , \widehat{X}_{i} , \dots , X_{+1}) \;\;\; + \\ \\
 & \sum_{i<j} (-1)^{i+j} \omega([X_{i},X_{j}], \dots ,\widehat{X}_i,
\dots ,\widehat{X}_j, \dots ,X_{q+1})
\end{array}
\vspace{1ex}
$$
v\'{e}rifie \( d_{\cal F}^{2} = 0 \) car la courbure de $\nabla$ est nulle.
La cohomologie
\( H^{\ast}({\cal F};\nu^{\ast}{\cal F}) \)
du complexe
\( (\Omega^{\ast}({\cal F};\nu^{\ast}{\cal F}),d_{\cal F}) \)
est la
{\bf cohomologie feuillet\'{e}e de $(M,{\cal F})$ \`{a}
valeurs dans
\( \nu^{\ast}{\cal F} \)}.
Alors, on a:
$$
E_{1}^{1,q}({\cal F}) \cong H^{q}({\cal F};\nu^{\ast}{\cal F})
$$

\subsection{Formes feuillet\'{e}es symplectiques}

\hspace{1em}
Une forme feuillet\'ee est un \'el\'ement d'un quotient,
mais on v\'erifie ais\'ement que ses
puissances ext\'erieures et son \'evaluation sur les vecteurs tangents
\`{a} $\cal F$ sont bien d\'{e}finies.
Une 2-forme feuillet\'{e}e
\( \sigma \in \Omega^{2}({\cal F}) \)
est dite
{\bf symplectique} si

\noindent
i) \( d_{\cal F} \sigma = 0 \)

\noindent
ii) \( \; \bigwedge^{k}\sigma \; \)
est non nulle en tout point de $M$ o\`{u} $dim{\cal F} = 2k$.

\noindent
Une telle forme feuillet\'{e}e d\'{e}finit une structure de Poisson
$\Lambda$ sur $M$ dont le feuilletage caract\'{e}ristique est $\cal F$.

A la structure de Poisson
\( \Lambda = ({\cal F},\sigma) \),
on associe les deux \'{e}l\'{e}ments de la suite spectrale de
Leray-Serre d\'{e}finis par
\vspace{1ex}

\noindent
i)
\( [\Lambda] = [\omega] \in E_{1}^{0,2}({\cal F}) = H^{2}({\cal F}) \)
\vspace{1ex}

\noindent
ii)
\( d_{1}[\Lambda] = [d_{1,0}\omega] \in E_{1}^{1,2}({\cal F}) =
H^{2}({\cal F};\nu^{\ast}{\cal F}) \)
\vspace{1ex}

\noindent
o\`{u} $\omega$ est un repr\'{e}sentant pur de type $(0,2)$
de $\sigma$ et $d_{1}$ est la diff\'{e}rentielle de la suite spectrale
induite par $d_{1,0}$. Ce sont des invariants essentiels de $\Lambda$.
\vspace{1ex}

\subsection{Alg\'{e}bro\"{\i}des de Lie}

\hspace{1em}
Une structure de Poisson $\Lambda=({\cal F},\sigma)$ sur une
vari\'{e}t\'{e} $M$
d\'{e}finit une structure d'{\bf alg\'{e}bro\"{\i}de de Lie} (\cite{Pr1})
sur le fibr\'{e} cotangent $T^{\ast}M$
(\cite{CDW}, \cite{Ka}) donn\'{e}e par:
\vspace{1ex}

\noindent
i) le morphisme de fibr\'{e}s vectoriels
\( \Lambda^{\#}\colon T^{\ast}M \rightarrow TM \)
qui, \`{a} toute 1-forme $\mu$, associe le champ
\( \Lambda^{\#}\mu \)
tangent \`{a} $\cal F$ tel que
\( i_{\textstyle \Lambda^{\#}\mu} \sigma = - \overline{\mu} \),
o\`{u} $\overline{\mu}$ est la classe feuillet\'{e}e de $\mu$;
\vspace{1ex}

\noindent
ii) le crochet de Lie
$$
\{\mu_{1},\mu_{2}\}
= i_{\textstyle \Lambda^{\#}\mu_{1}}d\mu_2 -
i_{\textstyle \Lambda^{\#}\mu_{2}}d\mu_{1} +
d\Lambda(\mu_{1},\mu_{2})
$$
pour lequel
\( \Lambda^{\#}\colon \Omega^{1}(M) \rightarrow {\fraktur X}(M) \)
est un morphisme d'alg\`{e}bres de Lie qui v\'erifie
$ \; \{\mu_{1},f\mu_{2}\} = f\{\mu_{1},\mu_{2}\} +
(L_{\textstyle \Lambda^{\#}\mu_{1}}f)\mu_{2} \; $
pour toute fonction
\( f \in {\cal C}^{\infty}(M) \).
\vspace{3ex}

Le noyau du morphisme $\Lambda^{\#}$ est le
fibr\'{e} conormal \( \nu^{\ast}{\cal F} \). C'est un
{\bf alg\'{e}bro\"{\i}de de Lie vectoriel}, i.e. un fibr\'{e}
en alg\`{e}bres de Lie ab\'{e}liennes.
De fa\c{c}on pr\'{e}cise, on a une {\bf extension} d'alg\'{e}bro\"{\i}des
de Lie

\begin{equation}       \label{exta}
0 \rightarrow
\nu^{\ast}{\cal F}
\longrightarrow
T^{\ast}M
\stackrel{\Lambda^{\#}}{\longrightarrow}
T{\cal F}
\rightarrow 0
\end{equation}

\noindent
L'alg\'ebro\"{\i}de de Lie $T{\cal F}$ agit sur le noyau
$\nu^{\ast}{\cal F}$ \`a l'aide
de la connexion partielle de Bott.
Pour tout couple $X_1$ et $X_2$ de champs tangents \`a
$\cal F$, on consid\`ere le couple $\mu_1= - i_{X_1} \omega$ et
$\mu_2 = - i_{X_2} \omega$
de 1-formes pures de type $(0,1)$ pour le choix d'un
suppl\'ementaire de $T{\cal F}$.
La 2-forme feuillet\'{e}e ferm\'{e}e
$\Omega \in \Omega^2({\cal F};\nu^\ast{\cal F})$
d\'efinie par
$$
\Omega(X_{1},X_{2}) = - \{\mu_{1},\mu_{2} \}_{1,0} =
i_{\textstyle X_{2}}i_{\textstyle X_{1}}d_{1,0} \omega
$$
repr\'esente une classe de $H^2({\cal F};\nu^\ast{\cal F})$ qui
caract\'erise l'extension:
\vspace{1ex}

\begin{prop}[\cite{DH}]         \label{alg}
L'alg\'{e}bro\"{\i}de de Lie $T^{\ast}M$ de $(M,\Lambda)$ est
l'extension de $T{\cal F}$ par $\nu^{\ast}{\cal F}$
(relative \`{a} la connexion partielle de Bott) de classe
\( d_{1}[\Lambda] \). $\Box$
\end{prop}

\subsection{Int\'{e}gration de Poisson} \label{intP}

\hspace{1em}
Soit $(M_0,{\cal F}_0)$ une vari\'et\'e feuillet\'ee (o\`u l'on modifie
les notations
pr\'ec\'edentes par l'adjonction d'un indice 0).
Le {\bf groupo\"{\i}de d'homotopie}
$\Pi_{1}({\cal F}_{0})$ (\cite{Pr1}) est le quotient de l'espace des chemins
contenus dans les feuilles de ${\cal F}_{0}$ (muni de la topologie
compact-ouvert $C^{\infty}$) par la relation d'homotopie dans les
feuilles de ${\cal F}_{0}$. C'est un groupo\"{\i}de de Lie (voir \cite{P})
dont l'espace total $M$ est s\'{e}par\'{e} si, et seulement si, tout cycle
\'{e}vanouissant de ${\cal F}_{0}$ est trivial (voir \cite{DH}). Les
projections source $\alpha_{0}$ et but $\beta_{0}$ d\'{e}finissent un m\^{e}me
feuilletage image r\'{e}ciproque
\( {\cal F} = \alpha_{0}^{\ast}{\cal F}_{0}
=  \beta_{0}^{\ast}{\cal F}_{0} \)
dont les feuilles sont les groupo\"{\i}des d'homotopie des
feuilles de ${\cal F}_{0}$.
Ce feuilletage est invariant par l'involution $\iota_{0}$ de
\( \Pi_{1}({\cal F}_{0}) \)
et sa trace sur $M_{0}$ est \'egale \`a ${\cal F}_{0}$.
\vspace{2ex}

\noindent
i) Les applications $\alpha_{0}$, $\beta_{0}$, $\iota_{0}$ induisent
des morphismes des complexes de formes relatives et de formes feuillet\'{e}es.
\vspace{1ex}

\noindent
ii) On fixe une d\'ecomposition
\( TM_0 = N{\cal F}_0 \oplus T{\cal F}_0 \)
et l'on note $p$ la projection
sur $T{\cal F}_0$. Le noyau de l'application
$$
(p \times p) {\scriptstyle \circ} (\beta_{0},\alpha_{0})_{\ast}\colon TM
\rightarrow TM_{0} \oplus TM_{0} \rightarrow
T{\cal F}_{0} \oplus T{\cal F}_{0}
$$
est un suppl\'{e}mentaire de $T{\cal F}$. On dira que ces
d\'{e}compositions de $TM$ et $TM_{0}$ sont {\bf adapt\'{e}es l'une \`{a}
l'autre}.
Dans des d\'{e}compositions adapt\'{e}es, les applications
$\alpha_{0}$, $\beta_{0}$, $\iota_{0}$ induisent des morphismes des
complexes de formes pures de type $(p,q)$.
\vspace{1ex}

\noindent
iii) Soit
\( \Omega^{\ast}({\cal F},{\cal F}_{0}) \)
le complexe des {\bf formes feuillet\'{e}es relatives} qui s'annulent sur
$M_0$. En utilisant des d\'{e}compositions adapt\'{e}es, on obtient
un isomorphisme de
\( \Omega^{0,q}(M,M_{0}) \)
dans
\( \Omega^{q}({\cal F},{\cal F}_{0}) \)
et une
{\bf suite spectrale de Leray-Serre relative}
$$
E_{2}^{p,q}({\cal F},{\cal F}_{0}) \Rightarrow H^{p+q}(M,M_{0})
$$

En outre, toute structure de Poisson
\( \Lambda_{0} = ({\cal F}_{0},\sigma_{0}) \)
sur $M_{0}$ se rel\`{e}ve en une structure de Poisson $\Lambda$ sur
$M$ qui en fait un {\bf groupo\"{\i}de de Poisson} au sens de \cite{W3}
(voir \cite{DH}). Cette structure de Poisson est d\'{e}termin\'{e}e
par le feuilletage $\cal F$ et la forme feuillet\'{e} relative
\( \sigma = \alpha_{0}^{\ast}\sigma_{0} - \beta_{0}^{\ast}\sigma_{0} \).

\setcounter{equation}{0}
\setcounter{th}{0}

\section{P\'eriodes sph\'eriques et groupo\"{\i}de structural}

\hspace{1em}
Le probl\`{e}me de l'int\'{e}gration symplectique de $(M,\Lambda)$ consiste
\`{a} construire un groupo\"{\i}de symplectique dont l'alg\'{e}bro\"{\i}de
de Lie est isomorphe \`{a} l'extension $T^{\ast}M$ de $T{\cal F}$ par
\( \nu^{\ast}{\cal F} \).
Ces deux derniers alg\'{e}bro\"{\i}des de Lie sont int\'{e}grables.
Le fibr\'{e}
conormal $\nu^{\ast}{\cal F}$ est un groupo\"{\i}de de Lie (plus
pr\'ecisement un fibr\'{e} en groupes) qui est visiblement l'int\'{e}gration
de l'alg\'{e}bro\"{\i}de de Lie $\nu^{\ast}{\cal F}$. En outre, le
groupo\"{\i}de d'homotopie $\Pi_{1}({\cal F})$ r\'{e}alise l'int\'{e}gration
de $T{\cal F}$. L'int\'{e}grabilit\'{e} de $(M,\Lambda)$
sera caract\'{e}ris\'{e}e en termes de {\bf p\'{e}riodes sph\'{e}riques}
de la classe $d_{1}[\Lambda]$ de l'extension. Leur construction
sera le premier but de cette section. On s'int\'eressera ensuite
aux structures diff\'erentiable et feuillet\'ee du groupo\"{\i}de
quotient $\cal G$ de $\nu^{\ast}{\cal F}$ par ces p\'eriodes sph\'eriques.
C'est dans la construction de ce groupo\"{\i}de que se retrouveront toutes
les difficult\'es et toutes les obstructions \`a l'int\'egrabilit\'e.
En parti\-culier, on prouvera le th\'eor\`eme~\ref{th1} au \S\ref{groupat}.
\vspace{1ex}

\subsection{Le groupo\"{\i}de d'homotopie $\Pi_{2}({\cal F})$}

\hspace{1em}
Soit $\cal F$ un feuilletage de dimension n et de codimension m
sur une vari\'{e}t\'{e} $M$. Soit
\( {\cal C}^{\infty}(S^{2},{\cal F}) \)
l'espace des applications diff\'{e}rentiables de $S^{2}$ dans les feuilles
de $\cal F$ muni de la topologie compact-ouvert
${\cal C}^{\infty}$. Si $N$ d\'{e}signe le p\^{o}le
nord de $S^{2}$, la projection
\( \hat{p}\colon
s \in {\cal C}^{\infty}(S^{2},{\cal F})
\longmapsto
s({\scriptstyle N}) \in M \)
est continue et ouverte.
Soit $\Pi_2({\cal F})$ le quotient de
\( {\cal C}^{\infty}(S^{2},{\cal F}) \)
par la relation d'homotopie dans les feuilles de $\cal F$ relative \`{a} $N$
et
\( q\colon {\cal C}^{\infty}(S^{2},{\cal F}) \rightarrow \Pi_2(\cal F) \)
la projection quotient correspondante. La projection induite
\( p\colon \Pi_2({\cal F}) \rightarrow M \)
est aussi continue et ouverte; sa fibre en un point $x$ est le groupe
d'homotopie
\( \pi_2(L_{x},x) \)
de la feuille $L_x$ passant par $x$.
Bref,
\( \Pi_2({\cal F}) \)
est un groupo\"{\i}de (plus pr\'ecisement un fibr\'{e} en groupes)
topologique appel\'{e} le
{\bf groupo\"{\i}de d'homotopie de $\cal F$ d'ordre 2}.
\vspace{0.5ex}

\subsection{Structure diff\'{e}rentiable sur $\Pi_2({\cal F})$} \label{stdif}

\hspace{1em}
Soit
\( s\colon (S^{2},N) \rightarrow (M,x) \)
une application diff\'{e}rentiable dont l'image est contenue
dans la feuille de $\cal F$ passant par $x$; on dira que $s$ est
une {\bf sph\`ere tangente \`a $\cal F$}.
Soit $E$ le fibr\'{e} vectoriel image r\'{e}ciproque de
\( \nu^{\ast}{\cal F} \)
par $s$. A l'aide d'une m\'{e}trique riemannienne, on construit
une application diff\'{e}rentiable $D$ d'un voisinage de la section nulle de
$E$ dans $M$ prolongeant s et transverse \`{a} $\cal F$, qui induit un
diff\'eomorphisme de la fibre de $N$ sur une transversale $V$ passant par $x$.
D'apr\`{e}s le th\'{e}or\`{e}me de stabilit\'{e} globale de Reeb, le
feuilletage image r\'{e}ciproque $D^{\ast}{\cal F}$ est trivialis\'{e}
par la connexion de Bott.
Bref, $s$ se prolonge en une application diff\'{e}rentiable
\( D\colon S^{2} \times V \rightarrow M \)
transverse \`{a} $\cal F$ telle que le feuilletage image r\'{e}ciproque
${\cal H} = D^{\ast}{\cal F}$
est le feuilletage horizontal en sph\`eres sur
\( S^{2} \times V \).
On dira que $D$ est une {\bf d\'{e}formation transverse de $s$}.
\vspace{2ex}

Toute d\'{e}formation transverse $D$ peut \^{e}tre \'{e}paissie en un
{\bf tube de sph\`{e}res tangentes}.
Pour cela, soit
\( (U;x_1, \dots , x_n,y_1, \dots ,y_m) \)
une carte locale distingu\'ee pour laquelle la transversale
\( x_1 = \dots = x_n =0 \)
est contenue dans $V$.
Soient
\( \varphi^{1}_{t} , \dots , \varphi^{n}_{t} \)
les flots des champs tangents
\( \frac{\partial}{\partial x_1} , \dots , \frac{\partial}{\partial x_n} \).
L'application diff\'{e}rentiable
\( T\colon S^{2} \times U \rightarrow M \)
d\'{e}finie par
$$
T(z,x_{1}, \dots ,x_{n},y_{1}, \dots ,y_{m}) =
\varphi^{1}_{x_{1}} {\scriptstyle \circ} \dots {\scriptstyle \circ}
\varphi^{n}_{x_{n}}
(D(z,y_{1}, \dots , y_{n}))
$$
prolonge la d\'eformation transverse $D$.
On dira que $T$ est un {\bf tube de sph\`{e}res tangentes}.
Tout tube $T$ d\'{e}finit des applications continues

\begin{equation}
\raisebox{-4ex}{
\begin{picture}(150,55)(0,0)
\put(50,33){\vector(3,-1){50}}
\put(50,44){\vector(1,0){50}}
\put(122,33){\vector(0,-1){15}}
\put(34,41){$U$}
\put(108,41){${\cal C}^{\infty}(S^{2},{\cal F})$}             \label{tube}
\put(108,5){$\Pi_2({\cal F})$}
\put(72,14){$\tau$}
\put(129,25){$q$}
\put(72,48){$\hat{\tau}$}
\end{picture}}
\end{equation}

\noindent
o\`u $\tau$ est une section de la projection $p$.
Ces applications $\tau$ d\'efinissent un atlas et
donc une structure diff\'{e}rentiable sur
\( \Pi_2({\cal F}) \)
pour laquelle $p$ est un diff\'{e}omorphis\-me local. Bref,
on a le r\'{e}sultat suivant:

\begin{prop}
Le groupo\"{\i}de $\Pi_2({\cal F})$ est un groupo\"{\i}de de Lie
\'{e}tal\'{e} sur $M$. $\Box$
\end{prop}

\begin{exems}   \label{pi2fi}
{\em
(1) Soit $\Pi_2(M)$ le groupo\"{\i}de
d'homotopie d'ordre 2 d'une vari\'et\'e $M$.
Soient $V$ un voisinage contractile d'un point $x_0$ et
$\gamma$ un chemin dans $V$ d'extre\-mit\'e $x_0$.
L'isomorphisme induit
$\gamma_\#\colon \pi_2(M,x) \rightarrow \pi_2(M,x_0) $
ne d\'epend que de l'origine $x$ de $\gamma$.
Les hom\'eomorphismes
$$
[s] \in p^{-1}(V) \mapsto
(s({\scriptstyle N}),\gamma_\#([s])) \in V \times \pi_2(M,x_0)
$$
d\'efinissent une structure de vari\'et\'e sur $\Pi_2(M)$
qui en fait un fibr\'{e} localement trivial
de fibre $\pi_2(M,x_0)$ et de groupe structural $\pi_{1}(M,x_0)$.
Evidemment, ce fibr\'{e} est trivial si $M$ est {\bf 2-simple},
i.e. l'action de $\pi_{1}(M,x_0)$ sur $\pi_{2}(M,x_0)$ est triviale.
\vspace{1ex}

\noindent
(2) Soient
\( \pi\colon M \rightarrow B \)
un fibr\'{e} localement trivial dont la fibre $F$ est 2-simple et $\cal F$ le
feuilletage d\'efini par $\pi$. Soit $\Pi_2(\pi)$
le fibr\'e associ\'e de fibre $\pi_2(F)$.
Le groupo\"{\i}de d'homotopie $\Pi_2({\cal F})$ est le fibr\'e
image r\'{e}ciproque
de $\Pi_2(\pi)$ par $\pi$.
\vspace{1ex}

\noindent
(3) On a une situation analogue si
$\cal F$ est un feuilletage d\'{e}fini par une submersion surjective
\( \pi\colon M \rightarrow B \)
\`a fibres 2-simples. Le groupo\"{\i}de
d'homotopie $\Pi_2({\cal F})$ est trivial en restriction aux fibres de $\pi$
et donc d\'efinit par projection un groupo\"{\i}de de Lie $\Pi_2(\pi)$
\'{e}tal\'{e} sur $B$. On dira que $\Pi_2(\pi)$ est le
{\bf groupo\"{\i}de d'homotopie d'ordre 2 de $\pi$}.}
\end{exems}

Les fibres de
\( \alpha_{0}\colon \Pi_{1}({\cal F}) \rightarrow M \)
sont les rev\^{e}tements universels des feuilles de $\cal F$ de projection
$\beta_{0}$. Celle-ci induit donc un isomorphisme
\( (\beta_{0})_{\ast}\colon \Pi_2(\alpha_{0}) \rightarrow
\Pi_2(\cal F) \).
D'apr\`{e}s les exemples (1) et (2) ci-dessus,
on a le r\'{e}sultat suivant:

\begin{prop}                \label{pi2}
Le groupo\"{\i}de d'homotopie
$\Pi_2({\cal F})$ est localement trivial en res\-triction \`{a}
chaque feuille de $\cal F$.
Si $\Pi_{1}({\cal F})$ est en plus localement trivial, alors il en est de
m\^{e}me  pour $\Pi_2({\cal F})$. $\Box$
\end{prop}

\subsection{Int\'{e}gration et p\'{e}riodes}

\hspace{1em}
Soient $\Lambda=({\cal F},\sigma)$ une structure de Poisson et
$\Omega \in \Omega^2({\cal F};\nu^\ast{\cal F})$
un repr\'{e}sentant de la classe $d_{1}[\Lambda]$.
Soit  \( D\colon S^{2} \times V \rightarrow M \)
une d\'{e}formation transverse d'une sph\`{e}re tangente
\( s\colon (S^{2},N) \rightarrow (M,x) \).
On d\'esigne par $\cal H$ le feuilletage horizontal en sph\`eres sur
$S^{2} \times V$.
La 2-forme feuillet\'{e}e ferm\'{e}e
\( D^{\ast}\Omega \in \Omega^{2}({\cal H};\nu^{\ast}{\cal H}) \)
est repr\'esent\'ee par une forme
pure $\delta$ de type $(1,2)$ pour la trivialisation canonique de
$T^\ast(S^2 \times V)$. L'int\'egration
sur les fibres du fibr\'e trivial
$S^2 \times V$ au-dessus de $V$
d\'efinit une 1-forme
$$
\int_{\textstyle D} \: \Omega
= - \!\!\!\!\!\! \int_{\textstyle S^{2}} \: \delta
$$
sur $V$.
Les propri\'et\'es de $- \!\!\!\!\! \int$
impliquent que celle-ci est ind\'ependante de la tri\-vialisation de
$T^\ast(S^2 \times V)$.

De fa\c{c}on analogue, la 2-forme feuillet\'{e}e $\Omega$
s'int\`{e}gre sur un tube $T$ prolongeant $D$
en une 1-forme basique $\int_{T}\Omega$ sur l'ouvert distingu\'{e}
correspondant $U$ qui \'{e}tend la 1-forme $\int_{D}\Omega$ sur
la transversale $V$.
\vspace{1ex}

On consid\`ere le morphisme d'int\'egration
$$
I_{\textstyle \Omega}\colon \Pi_2({\cal F}) \rightarrow \nu^{\ast}{\cal F}
$$
qui, \`a la section
\( \tau\colon U  \rightarrow \Pi_2({\cal F}) \)
d\'efinie par un tube $T$ d'apr\`es le diagramme~(\ref{tube}),
associe la 1-forme basique $\int_T \Omega$ sur $U$.
D'apr\`es le th\'eor\`eme de Stokes, c'est un morphisme
de groupo\"{\i}des de Lie bien d\'efini
qui ne d\'epend que de la classe $d_1[\Lambda]$ de $\Omega$.
Son image
\( {\cal P} = {\cal P}(\Lambda) \)
est un sous-groupo\"{\i}de alg\'ebrique de
\( \nu^{\ast}{\cal F} \)
appel\'e le {\bf groupo\"{\i}de des p\'{e}riodes
sph\'{e}riques d\'{e}riv\'{e}es}
de $\Lambda$.
\vspace{1ex}

\begin{prop}  \label{sat}
Le groupo\"{\i}de d\'{e}riv\'{e} $\cal P$ est satur\'{e} pour la
connexion partielle de Bott sur
\( \nu^{\ast}{\cal F} \).
\end{prop}

\noindent
{\bf D\'{e}monstration}
Puisque les int\'{e}grales de $\Omega$ sur les tubes sont
des 1-formes basiques, le groupo\"{\i}de d\'{e}riv\'{e}
$\cal P$ est horizontal pour la connexion partielle
de Bott. D'apr\`{e}s la proposition~\ref{pi2}, $\Pi_2({\cal F})$ est
localement trivial en restriction aux feuilles de $\cal F$ et donc il
en est de m\^{e}me pour $\cal P$. D'o\`{u} la proposition. $\Box$

\subsection{Structure diff\'{e}rentiable sur $\cal P$}

\hspace{1em}
Si $\cal K$ d\'esigne le noyau du morphisme d'int\'egration, on a
une suite exacte de groupo\"{\i}des topologiques
$$
0 \rightarrow {\cal K} \longrightarrow
\Pi_2({\cal F})
\stackrel{I_{\Omega}}{\longrightarrow} {\cal P} \rightarrow 0
$$
o\`u le groupo\"{\i}de d\'eriv\'e
$ {\cal P} $
est muni de la topologie quotient.

L'exemple~\ref{exobs}.1 sugg\`ere la d\'efinition suivante:
une d\'{e}formation transverse $D$ d'une sph\`ere tangente
\( s\colon (S^2,N) \rightarrow (M,x) \)
est {\bf symplectiquement
\'{e}vanouis\-sante} si
\( \int_{D}\Omega (x) = 0 \);
une telle d\'{e}formation est {\bf triviale} si
\( \int_{D} \Omega = 0 \).
\vspace{1ex}

\begin{prop}            \label{defsym}
Les conditions suivantes sont \'{e}quivalentes:

\noindent
i) le groupo\"{\i}de $\cal P$ est
un sous-groupo\"{\i}de de Lie de
\( \nu^{\ast}{\cal F} \)
\'{e}tal\'{e} sur $M$;

\noindent
ii) le noyau $\cal K$ est un sous-groupo\"{\i}de de Lie
plong\'{e} de $\Pi_2({\cal F})$;

\noindent
iii) toute d\'{e}formation symplectiquement \'{e}vanouissante est triviale.
\end{prop}

\noindent
{\bf D\'{e}monstration}
On montre tout d'abord que les conditions (i) et (ii) sont \'{e}quivalentes
(cf. \cite{HM} et \cite{Pr2}).
Si $\cal P$ est un sous-groupo\"{\i}de de Lie de
\( \nu^{\ast}{\cal F} \)
\'{e}tal\'{e} sur $M$, le morphisme d'int\'{e}gration
\( I_{\Omega}\colon \Pi_2({\cal F}) \rightarrow {\cal P} \)
est un diff\'{e}omorphisme local. Son noyau
${\cal K} = I_\Omega^{-1} (M)$ est donc un
sous-groupo\"{\i}de de Lie plong\'{e} de
\( \Pi_2({\cal F}) \).

R\'{e}ciproquement, soit $\cal R$ la relation d'\'{e}quivalence sur
\( \Pi_2({\cal F}) \)
d\'{e}finie par l'action de $\cal K$,
c'est-\`{a}-dire l'image r\'{e}ciproque de $\cal K$
par l'application diff\'{e}rence
$$
([s_1],[s_2]) \in \Pi_2({\cal F}) \times_M \Pi_2({\cal F})
\longmapsto
[s_2] - [s_1] \in \Pi_2({\cal F})
$$
o\`u
\( \Pi_2({\cal F}) \times_M \Pi_2({\cal F}) \)
est le produit fibr\'e
\( \{ ([s_1],[s_2]) \colon s_1({\scriptstyle N}) = s_2({\scriptstyle N}) \} \).
La condition (ii) implique que $\cal R$ est une sous-vari\'{e}t\'{e}
plong\'{e}e de
\( \Pi_2({\cal F}) \times_M \Pi_2({\cal F}) \).
De plus, la projection sur le premier facteur se restreint
en une submersion surjective de $\cal R$ sur $\Pi_2({\cal F})$.
D'apr\`{e}s le crit\`{e}re de Godement (voir \cite{S}), le quotient
\( {\cal P} =\Pi_2({\cal F}) / {\cal K} \)
est muni d'une structure de vari\'{e}t\'{e} (qui en fait
un groupo\"{\i}de de Lie) pour laquelle la projection $I_\Omega$ est une
submersion surjective. L'inclusion de $\cal P$
dans $\nu^\ast{\cal F}$ est un morphisme injectif de groupo\"{\i}des de Lie.
Puisque $\Pi_2({\cal F})$ s'\'etale sur $M$, il en est de m\^eme pour
$\cal P$ qui devient ainsi immerg\'e dans $\nu^\ast{\cal F}$.
Bref, $\cal P$ est un sous-groupo\"{\i}de de Lie de
\( \nu^{\ast}{\cal F} \)
\'etal\'e sur $M$.

Pour montrer l'\'{e}quivalence entre (ii) et (iii), on consid\`{e}re une
d\'{e}formation transverse $D$ d'une sph\`{e}re tangente
$s$. Un tube $T$ qui prolonge $D$ d\'{e}finit un voisinage
\( W = \tau (U) \)
de $[s]$ diff\'{e}omorphe \`a l'ouvert
distingu\'{e} $U$ associ\'e \`a $T$.
L'\'{e}quivalence est une cons\'{e}quence
imm\'{e}diate des remarques suivantes:

\noindent
1) la classe
\( [s] \in {\cal K} \)
si et seulement si la d\'{e}formation $D$ est symplectiquement
\'{e}vanouissante;

\noindent
2) le voisinage $W$ de $[s]$ est contenu dans $\cal K$ si et seulement si
l'int\'{e}grale $\int_{T}\Omega$ est nulle. Or cette 1-forme basique est
nulle si et seulement si
\( \int_{D} \Omega = 0 \). $\Box$
\vspace{1ex}

\subsection{Groupo\"{\i}de structural: structure diff\'erentiable}
\label{groupat}

\hspace{1em}
Soit $\cal P$ le groupo\"{\i}de des p\'eriodes sph\'eriques
d\'{e}riv\'{e}es de $\Lambda$.
On consid\`{e}re la suite exacte de groupo\"{\i}des topologiques

\begin{equation}        \label{seg}
0 \rightarrow {\cal P} \longrightarrow
\nu^{\ast}{\cal F}
\stackrel{q}{\longrightarrow} {\cal G} \rightarrow 0
\end{equation}

\noindent
o\`{u} $\cal G$ est le groupo\"{\i}de quotient
\( \nu^{\ast}{\cal F} / {\cal P} \). On dira que $\cal G$ est
le {\bf groupo\"{\i}de
structural} de la vari\'et\'e de Poisson $(M,\Lambda)$.
En proc\'edant comme dans la preuve de
la proposition~\ref{defsym}, on d\'emontre le th\'eor\`eme~\ref{th1},
\`a savoir que
$\cal G$ est un groupo\"{\i}de de Lie si et seulement si
le groupo\"{\i}de d\'eriv\'e $\cal P$ est un sous-groupo\"{\i}de
de Lie plong\'e de $\nu^\ast{\cal F}$ \'etal\'e sur $M$.
\vspace{1ex}

\subsection{Groupo\"{\i}de structural: structure feuillet\'ee}

\hspace{1em}
On supposera d\'esormais que $\cal G$ est
un groupo\"{\i}de  de Lie. Soient $L$ la forme de Liouville sur le fibr\'e
conormal  \( p\colon \nu^{\ast}{\cal F} \rightarrow M \)
et $\cal S$ le feuilletage image r\'{e}ciproque
de $\cal F$ par $p$. La 1-forme feuillet\'{e}e
\( \Phi \in \Omega^{1}({\cal S};\nu^{\ast}{\cal S}) \)
d\'{e}finie par
$$
\Phi (X) = i_{\textstyle X} dL
$$
sera appel\'ee la {\bf forme de Maurer-Cartan} de
\( \nu^{\ast}{\cal F} \).
Son noyau est le sous-fibr\'{e} horizontal de la connexion partielle de Bott
d'apr\`es~(\ref{rel}). Cette forme feuillet\'ee v\'erifie les propri\'et\'es
suivantes que l'on d\'{e}duit ais\'{e}ment de l'\'{e}criture locale de la
forme de Liouville $L$:
\vspace{1ex}

\noindent
1)
\( \:\: \Phi(\mu^{\ast}) = p^{\ast}\mu \:\: \)
pour tout champ invariant \`{a} gauche $\mu^{\ast}$
associ\'{e} \`{a} une section $\mu$ de l'alg\'{e}bro\"{\i}de de Lie
\( \nu^{\ast}{\cal F} \).
\vspace{1ex}

\noindent
2) $\Phi$ est {\bf invariante \`{a} gauche},
i.e. pour toute 1-forme basique $\mu$, on a
$$
(L_{\textstyle \mu})^{\ast} \Phi = \Phi
$$
o\`u $L_{\textstyle \mu}$ est la {\bf translation \`a gauche}
pour la structure de groupo\"{\i}de sur $\nu^\ast{\cal F}$;
\vspace{1ex}

\noindent
3)
\( \:\: d_{\cal S} \Phi = 0 \).
\vspace{1ex}

\noindent
La forme de Maurer-Cartan v\'erifie en fait la propri\'et\'e
d'invariance suivante: pour toute
section $\mu$ de $\nu^\ast{\cal F}$, on a

\begin{equation}    \label{Lmu}
(L_{\textstyle \mu})^{\ast} \Phi  - \Phi = p^\ast d_{\cal F} \mu
\end{equation}

\noindent
et donc

\begin{equation}   \label{prouni}
 \mu^\ast \Phi  = d_{\cal F} \mu
\end{equation}

Puisque $\cal P$ est horizontal pour la connexion partielle de
Bott, tous les objets pr\'ec\'edents sont projetables par
\( q\colon \nu^{\ast}{\cal F} \rightarrow {\cal G} \)
et d\'efinissent des objets analogues sur $\cal G$.
On dira que le groupo\"{\i}de de Lie
ab\'elien $\cal G$ est {\bf $\cal F$-feuillet\'e}.
\vspace{1ex}

\subsection{Compl\'etion transverse} \label{comp}

\hspace{1em}
On d\'esigne encore par $\cal S$ le feuilletage image r\'eciproque de
$\cal F$ par la projection
\( p\colon {\cal G} \rightarrow M \).
Si le groupo\"{\i}de structural $\cal G$ est s\'epar\'e,
on peut toujours construire un
suppl\'ementaire $N{\cal S}$ de
$T{\cal S}$ \`a l'aide d'une m\'etrique riemannienne.
Si $\cal G$ n'est pas s\'epar\'e,
il n'existe pas en g\'en\'eral de suppl\'ementaire de $T{\cal S}$.
Cela justifie la d\'efinition suivante: on dira que la vari\'et\'e de Poisson
$(M,\Lambda)$ est {\bf transversalement compl\`ete} si le
groupo\"{\i}de structural $\cal G$ (resp. le groupo\"{\i}de d\'eriv\'e
$\cal P$) poss\`ede un suppl\'ementaire
$N{\cal S}$ de $T{\cal S}$ (resp. invariant par l'action de $\cal P$).

\begin{exems}
{\em
(1) Soit $\cal F$ le feuilletage de
\( M = S^{2} \times {\Bbb R} - \{ (N,0) \} \)
d\'efini par l'\'equation $dt=0$.
En multipliant la forme volume canonique $v_0$ de $S^2$ par une
fonction positive convenable,
on peut obtenir une forme volume des feuilles $\omega$
dont la fonction d'aire
\( - \!\!\!\!\! \int_{S^2} \omega \;
\in \; {\cal C}^\infty({\Bbb R} - \{0\}) \)
tend vers $+\infty$ quand $t$ tend vers $0$.
Soit $\sigma$ la forme feuillet\'ee repr\'esent\'ee par $\omega$.
Le groupo\"{\i}de d\'eriv\'e de la structure de Poisson
$\Lambda=({\cal F},\sigma)$ est engendr\'e par la p\'eriode
$$
- \!\!\!\!\! \int_{S^2} d_{0,1}\omega
= d ( - \!\!\!\!\! \int_{S^2} \omega ) \;
\in \; \Omega^1({\Bbb R} - \{0\})
$$
C'est un sous-groupo\"{\i}de de Lie ferm\'e de $\nu^\ast{\cal F}$.
\vspace{1ex}

\noindent
(2) Si l'on remplace la 2-forme $\omega$ de l'exemple (1) par la 2-forme
$e^{\textstyle t}v_{0}$, le groupo\"{\i}de
d\'eriv\'e n'est plus ferm\'e, mais le suppl\'ementaire
de $T{\cal S}$ engendr\'e par le champ
\( \frac{\partial}{\partial t} + r \frac{\partial}{\partial r} \)
est invariant par l'action de $\cal P$.
\vspace{1ex}

Ces deux vari\'et\'es de Poisson sont int\'egrables d'apr\`es le
th\'eor\`eme~\ref{th2}. On exhibe enfin un exemple de vari\'et\'e de Poisson
qui n'est pas transversalement compl\`ete:
\vspace{1ex}

\noindent
(3) Soit $M$ le compl\'ementaire dans $S^2 \times S^2 \times {\Bbb R}$
de la r\'eunion des
deux ensembles
\( S^{2} \times \{N\} \times \rbrack -\infty,0 \rbrack \)
et
\( \{N\} \times S^{2} \times \lbrack 0,+\infty \lbrack \).
Le feuilletage en produit de sph\`eres
induit un feuilletage $\cal F$ d\'efini par une submersion
\( f\colon M \rightarrow {\Bbb R} \).
Soient $p_1$ et $p_2$ les projections de $S^2 \times S^2$ sur chacun
des facteurs. La forme
feuillet\'{e}e symplectique $\sigma$ repr\'{e}sent\'{e}e par
\( e^{-t} \; p_{1}^{\ast} v_{0} \; +\; e^t \;p_{2}^{\ast} v_{0} \)
d\'{e}finit une structure de
Poisson $\Lambda$ sur $M$.
Le groupo\"{\i}de d\'{e}riv\'{e} $\cal P$ est engendr\'{e}
par la 1-forme $e^{-t} dt$ au-dessus de
$M_- = f^{-1}(\rbrack -\infty,0 \lbrack)$
et par la 1-forme $me^t dt$ au-dessus de
$M_+ = f^{-1}(\rbrack 0,+\infty \lbrack)$.
La forme conormale dt en tout point de la feuille
$f^{-1}(0)$ est simultan\'eement adh\'erente aux nappes de
p\'eriodes d\'efinies par $e^{-t} dt$ au-dessus de $M_-$ et $e^t dt$
au-dessus de $M_+$. L'invariance d'un suppl\'ementaire $N{\cal S}$
par l'action de ces deux nappes entra\^{\i}nerait que le champ vertical
$\frac{\partial}{\partial r}$ est normal en restriction \`a
$f^{-1}(0)$ ce qui n'est \'evidemment pas possible.}
\end{exems}

\subsection{Groupo\"{\i}de structural de l'int\'egration de Poisson}

\hspace{1em}
Soient $\Lambda_0=({\cal F}_0,\sigma_0)$ une structure de Poisson
sur une vari\'et\'e
$M_0$ et $\Lambda=({\cal F},\sigma)$ la structure de Poisson relev\'ee sur
\( M = \Pi_1({\cal F}_0) \).
Le fibr\'{e} conormal
\( \nu^{\ast}{\cal F} \)
est le fibr\'{e} vectoriel image r\'{e}ciproque du fibr\'{e}
conormal
\( \nu^{\ast}{\cal F}_{0} \)
par $\alpha_{0}$ et $\beta_{0}$.
On se propose de montrer qu'il en est de m\^eme pour les
groupo\"{\i}des structuraux
$\cal G$ et ${\cal G}_0$ de $\Lambda$ et $\Lambda_0$.
Pour cela, il suffit de d\'emontrer le lemme suivant:
\vspace{1ex}

\begin{lem}    \label{perhaut}
Le groupo\"{\i}de d\'{e}riv\'{e} $\cal P$ de $\Lambda$ est
l'image r\'{e}ciproque du groupo\"{\i}de d\'{e}riv\'{e}
${\cal P}_{0}$  de $\Lambda_{0}$ par $\alpha_{0}$ et $\beta_{0}$.
\end{lem}

\noindent
{\bf D\'{e}monstration} Si $\Omega_0$ est un repr\'{e}sentant de la classe
$d_1[\Lambda_0]$, alors la classe $d_1[\Lambda]$ est repr\'{e}sent\'{e}e par
\( \Omega = \alpha_{0}^{\ast} \Omega_{0} - \beta_{0}^{\ast} \Omega_{0} \).
Pour toute sph\`{e}re $s$ tangente \`a la feuille de $\cal F$
passant par $u \in M_{0}$, la p\'{e}riode

\begin{equation}  \label{per1}
\int_{\textstyle s} \; \Omega =
\alpha_{0}^{\ast} ( \int_{\textstyle \alpha_{0} {\textstyle \circ} s}
\; \Omega_{0} ) -
\beta_{0}^{\ast} ( \int_{\textstyle \beta_{0} {\textstyle \circ} s }
\; \Omega_{0} )
\end{equation}

\noindent
D'autre part, toute sph\`{e}re $s_{0}$ tangente \`a la feuille $L_0$ de
${\cal F}_{0}$ passant par $u$ se rel\`{e}ve par $\beta_{0}$ en une sph\`{e}re
$s$ tangente au rev\^{e}tement universel $\alpha_{0}^{-1}(u)$ de $L_0$.
Leurs p\'{e}riodes sph\'{e}riques sont reli\'{e}es par

\begin{equation}   \label{per2}
\int_{\textstyle s} \; \Omega = - \beta_{0}^{\ast}
( \int_{\textstyle s_{0}} \; \Omega_{0} )
\end{equation}

\noindent
D'apr\`{e}s (\ref{per1}) et (\ref{per2}), ${\cal P}_{0}$ est la
restriction de $\cal P$ \`{a} $M_{0}$.
Par ailleurs, la connexion partielle de Bott trivialise $\cal P$ en
restriction aux fibres de $\alpha_{0}$ et de $\beta_0$ car celles-ci sont
simplement connexes. D'o\`{u} le r\'{e}sultat. $\Box$
\vspace{1ex}

\begin{prop}    \label{ghaut}
Le groupo\"{\i}de structural
$ p\colon {\cal G} \to  M $ de $\Lambda$ est l'image r\'{e}ciproque du
groupo\"{\i}de
structural $p_0\colon {\cal G}_{0} \to M_0$ de $\Lambda_0$ par $\alpha_0$ et
$\beta_0$. En outre, si $\cal S$ et ${\cal S}_0$ sont les feuilletages image
r\'{e}ciproque de $\cal F$ et ${\cal F}_0$
par $p$ et $p_0$, alors la forme de Maurer-Cartan
\( \Phi \in \Omega^{1}({\cal S};\nu^{\ast}{\cal S}) \)
de $\cal G$ est l'image r\'{e}ciproque de la forme
de Maurer-Cartan
\( \Phi_0 \in \Omega^{1}({\cal S}_0;\nu^{\ast}{\cal S}_0) \)
de ${\cal G}_0$ par $\alpha_0$ et $\beta_0$. $\Box$
\end{prop}

\setcounter{equation}{0}
\setcounter{th}{0}

\section{Int\'egration symplectique}

\hspace{1em}
Soit $\Lambda_{0}=({\cal F}_{0},\sigma_{0})$ une structure de
Poisson sur une vari\'{e}t\'{e} $M_{0}$. Dans cette section, on va d\'emontrer
le th\'eor\`eme~\ref{th2}. On rappelle son \'enonc\'e: {\em si
\vspace{2ex}

\noindent
$H_1$) tout cycle \'{e}vanouissant de ${\cal F}_0$ est trivial,
\vspace{1ex}

\noindent
$H_2$) le groupo\"{\i}de structural
${\cal G}_0$ est un
groupo\"{\i}de de Lie,
\vspace{1ex}

\noindent
$H_3$) la vari\'et\'e de Poisson $(M_0,\Lambda_0)$ est transversalement
compl\`ete,  \vspace{2ex}

\noindent
alors $(M_{0},\Lambda_{0})$ est int\'{e}grable.}
\vspace{3ex}

Soit $\Lambda = ({\cal F},\sigma)$ la structure de Poisson relev\'{e}e
sur $M=\Pi_1({\cal F}_0)$. Int\'egrer $(M_{0},\Lambda_{0})$ va consister \`a
compl\'{e}ter l'int\'{e}gration de Poisson
$(M,\Lambda)$ en un diagramme

\begin{center}
\begin{picture}(260,53)(0,0)
\put(62,33){\vector(2,-1){35}}
\put(59,44){\vector(1,0){44}}
\put(143,44){\vector(1,0){44}}
\put(184,33){\vector(-2,-1){35}}
\put(192,33){\vector(-2,-1){35}}
\put(119,33){\vector(0,-1){20}}
\put(125,33){\vector(0,-1){20}}
\put(15,41){$({\cal G}_{0},\eta_{0})$}
\put(195,41){$(M,\Lambda)$}
\put(109,41){$(\Gamma,\eta)$}
\put(104,0){$(M_{0},\Lambda_{0})$}
\put(67,18){$p_{0}$}
\put(109,24){$\alpha$}
\put(129,24){$\beta$}
\put(152,28){$\alpha_{0}$}
\put(174,16){$\beta_{0}$}
\put(78,48){$i$}
\put(160,48){$\pi$}
\end{picture}
\end{center}
\vspace{0.5ex}

\noindent
o\`u $(\Gamma,\eta)$ est un groupo\"{\i}de symplectique,
$({\cal G}_{0},\eta_{0})$
est un sous-groupo\"{\i}de de Poisson {\bf pr\'esymplectique} (voir \cite{DH})
et $i$ et $\pi$ sont des morphismes de Poisson.

De fa\c{c}on pr\'{e}cise, on cherche \`{a} construire une {\bf extension}
$\Gamma$ du groupo\"{\i}de d'homotopie $M$
par le groupo\"{\i}de structural ${\cal G}_{0}$ qui
{\bf int\`{e}gre} l'extension d'alg\'{e}\-bro\"{\i}des de Lie~(\ref{exta}).
L'hypoth\`ese ($H_2$)
fournira donc le noyau de l'extension. On remarque que le groupo\"{\i}de
structural $\cal G$ de $(M,\Lambda)$ sera de m\^eme un groupo\"{\i}de de Lie
ab\'{e}lien
$\cal F$-feuillet\'{e} d'apr\`es la proposition~\ref{ghaut}.
La d\'emonstration se fera en trois \'etapes:
\vspace{2ex}

\noindent
1) {\em Int\'egration cohomologique}:
on montrera que la classe
\( d_{1}[\Lambda] \in H^{2}({\cal F},{\cal F}_{0};\nu^{\ast}{\cal F}) \)
{\bf s'int\`egre} en une classe
\( \nu \in H^{2}(M,M_{0};{\cal P}) \); celle-ci sera la classe
d'un {\bf cocycle sur $M$ \`{a} valeurs
dans $\cal G$} qui est cohomologue \`{a} z\'{e}ro en restriction \`{a} $M_{0}$.
L'hypoth\`ese ($H_1$) (i.e. la s\'eparation de $M$) interviendra dans
cette \'etape.
\vspace{2ex}

\noindent
2) {\em Int\'egration diff\'erentiable}:
cette \'etape se d\'ecomposera en trois nivaux
\vspace{1ex}

\noindent
i) on construira \`a l'aide du cocycle un {\bf fibr\'{e} principal
\( \pi\colon \Gamma \rightarrow M \)
de groupo\"{\i}de structural $\cal G$}
qui sera trivial en restriction \`a $M_{0}$;
\vspace{1ex}

\noindent
ii) on construira une {\bf connexion $\cal F$-partielle} $\Theta$ sur le
$\cal G$-fibr\'e
principal $\Gamma$;
\vspace{1ex}

\noindent
iii) on montrera que la {\bf courbure} $\Omega$ de $\Theta$ repr\'{e}sente
la classe
$d_{1}[\Lambda]$ qui deviendra ainsi une vraie {\bf classe de Chern.}
\vspace{2ex}

\noindent
3) {\em Int\'egration symplectique}:
il y aura encore deux niveaux diff\'erents
\vspace{1ex}

\noindent
i) on construira une forme symplectique
$\eta$ telle que
\( \pi\colon (\Gamma,\eta) \rightarrow (M,\Lambda) \)
soit un morphisme de Poisson \`a fibres isotropes.
L'hypoth\`ese ($H_3$) n'interviendra que pour construire $\eta$ dans
le cas o\`u l'espace total
$\Gamma$ n'est pas s\'epar\'e.
\vspace{1ex}

\noindent
ii) on montrera que les projections
$\alpha$ et $\beta$ forment une {\bf paire duale compl\`ete} (au sens de
\cite{W1}) ce qui d\'efinira une structure de groupo\"{\i}de symplectique sur
$\Gamma$ d'apr\`es \cite{CDW}.
\vspace{1ex}

\subsection{Int\'{e}gration cohomologique}

\hspace{1em}
Un {\bf cocycle sur $M$ \`{a} valeurs dans $\cal G$} est un couple
\( (\{U_i\},\{\tau_{ij}\}) \)
form\'{e} d'un recouvrement ouvert
\( {\cal U} = \{ U_i \} \)
et de sections
$\tau_{ij}\colon U_i \cap U_j \rightarrow {\cal G}$
telles que
\vspace{1ex}

\noindent
i) \( \: \tau_{ii} = 0 \: \) sur $U_i$,
\vspace{1ex}

\noindent
ii)
\( \: \tau_{ij} - \tau_{ik} + \tau_{jk} = 0 \: \)
sur $U_i \cap U_j \cap U_k$.
\vspace{1ex}

\noindent
Un tel cocycle repr\'{e}sente une classe dans le groupe
\( H^{1}({\cal U};
\raisebox{-1.1ex}{$\stackrel{\textstyle {\cal G}}{\sim}$}) \)
de cohomologie de $\cal U$
\`{a} valeurs dans le faisceau
\( \raisebox{-1.1ex}{$\stackrel{\textstyle {\cal G}}{\sim}$} \)
des germes de sections de $\cal G$.
Cette classe d\'{e}finit par passage \`a la limite une classe
$\tau$ dans le groupe \( \check{H}^{1}(M;
\raisebox{-1.1ex}{$\stackrel{\textstyle {\cal G}}{\sim}$}) \)
de la cohomologie de \v{C}ech de $M$
\`{a} valeurs dans le faisceau
$\raisebox{-1.1ex}{$\stackrel{\textstyle {\cal G}}{\sim}$}$ (voir
\cite{G}).
De fa\c{c}on g\'en\'erale, cette cohomologie v\'erifie les propri\'et\'es
suivantes:
\vspace{2ex}

\noindent
i) Si les cycles \'evanouissants de ${\cal F}_0$ sont triviaux,
la vari\'{e}t\'{e} $M$ est  s\'{e}par\'{e}e et donc la cohomologie de \v{C}ech
\( \check{H}^{\ast}(M;
\raisebox{-1.1ex}{$\stackrel{\textstyle {\cal G}}{\sim}$}) \)
est isomorphe \`{a} la cohomologie
\( H^{\ast}(M;
\raisebox{-1.1ex}{$\stackrel{\textstyle {\cal G}}{\sim}$}) \)
de $M$ \`{a} valeurs dans le faisceau
$\raisebox{-1.1ex}{$\stackrel{\textstyle {\cal G}}{\sim}$}$ (voir \cite{G}).
\vspace{2ex}

\noindent
ii) Pour toute paire $(M,M_{0})$, on a une suite exacte longue
$$
\dots \rightarrow
H^{q}(M,M_{0};
\raisebox{-1.1ex}{$\stackrel{\textstyle {\cal G}}{\sim}$})
\longrightarrow
H^{q}(M;
\raisebox{-1.1ex}{$\stackrel{\textstyle {\cal G}}{\sim}$})
\stackrel{\varepsilon_{0}^{\ast}}{\longrightarrow}
H^{q}(M_{0};
\raisebox{-1.1ex}{$\stackrel{\textstyle {\cal G}_0}{\sim}$})
\rightarrow \dots
$$
et puisque les applications induites par $\alpha_{0}$ et $\beta_{0}$ sont des
sections de $\varepsilon_{0}^{\ast}$,
celle-ci se d\'ecompose en suites exactes courtes:

\begin{equation}  \label{cohrel}
0 \rightarrow
H^{q}(M,M_{0};
\raisebox{-1.1ex}{$\stackrel{\textstyle {\cal G}}{\sim}$})
\rightarrow
H^{q}(M;
\raisebox{-1.1ex}{$\stackrel{\textstyle {\cal G}}{\sim}$})
\rightarrow
H^{q}(M_{0};
\raisebox{-1.1ex}{$\stackrel{\textstyle {\cal G}_0}{\sim}$})
\rightarrow 0
\end{equation}
\vspace{0.2ex}

\noindent
iii) La suite exacte de groupo\"{\i}des de Lie~(\ref{seg})
induit une suite exacte de faisceaux
\( 0
\rightarrow
{\cal P}
\rightarrow
\raisebox{-1.1ex}{$\stackrel{\textstyle \nu}{\sim}$}^{\ast}{\cal F}
\stackrel{q}{\rightarrow}
\raisebox{-1.1ex}{$\stackrel{\textstyle {\cal G}}{\sim}$}
\rightarrow
0 \).
Celle-ci induit une suite exacte longue
$$
\dots \rightarrow
H^{q}(M,M_{0};{\cal P})
\rightarrow
H^{q}(M;M_0;
\raisebox{-1.1ex}{$\stackrel{\textstyle \nu}{\sim}$}^{\ast}{\cal F})
\rightarrow
H^{q}(M,M_{0};
\raisebox{-1.1ex}{$\stackrel{\textstyle \cal G}{\sim}$})
\rightarrow \dots
$$
dont le cobord
\( \delta\colon H^{q}(M,M_{0};
\raisebox{-1.1ex}{$\stackrel{\textstyle {\cal G}}{\sim}$})
\rightarrow
H^{q+1}(M,M_{0};{\cal P}) \)
est un isomorphisme car le faisceau
\( \raisebox{-1.1ex}{$\stackrel{\textstyle \nu}{\sim}$}^{\ast}{\cal F} \)
des germes de 1-formes relatives est mou (voir \cite{G}).
\vspace{2ex}

\noindent
iv) Soient $\phi^1$ le faisceau des germes de 1-formes basiques et $\cal Q$
le faisceau quotient $\phi^1 / {\cal P}$. La suite exacte de faisceaux
\( 0 \rightarrow {\cal P}  \stackrel{j}{\longrightarrow}
\phi^{1}
\stackrel{k}{\longrightarrow}
{\cal Q} \rightarrow 0 \)
induit une suite exacte longue

\begin{equation}         \label{coQ}
\dots \rightarrow
H^{q}(M,M_0;{\cal P})
\stackrel{j^{\ast}}{\longrightarrow}
H^{q}(M,M_0;\phi^{1})
\stackrel{k^{\ast}}{\longrightarrow}
H^{q}(M,M_0;{\cal Q})
\rightarrow \dots
\end{equation}

Soit $F_u$ la fibre de la submersion
\( \alpha_{0}\colon M \rightarrow M_{0} \)
en un point $u $ de $M_0$. Les inclusions des fibres $F_u$
dans $M$ induisent des morphismes de restriction qui rendent commutatif
le diagramme suivant:

\begin{equation}
\raisebox{-5ex}{
\begin{picture}(330,45)(0,0)
\put(84,43){\vector(1,0){40}}
\put(84,3){\vector(1,0){40}}
\put(10,40){$H^{2}(M,M_0;{\cal P})$}
\put(50,32){\vector(0,-1){15}}
\put(2,0){$\prod H^{2}(F_u;{\cal P} \!\! \mid_{\textstyle F_u})$}
\put(214,43){\vector(1,0){40}}
\put(214,3){\vector(1,0){40}}
\put(136,40){$H^{2}(M,M_0;\phi^{1})$}
\put(268,40){$H^{2}(M,M_0;{\cal Q})$}                 \label{diares}
\put(172,32){\vector(0,-1){15}}
\put(304,32){\vector(0,-1){15}}
\put(130,0){$\prod H^{2}(F_u;\phi^{1} \!\! \mid_{\textstyle F_u})$}
\put(260,0){$\prod H^{2}(F_u;{\cal Q} \!\! \mid_{\textstyle F_u})$}
\put(92,8){$\{j^{\ast}_u\}$}
\put(103,49){$j^{\ast}$}
\put(230,49){$k^{\ast}$}
\put(221,8){$\{k^{\ast}_u\}$}
\put(58,22){$\rho_{\cal P}$}
\put(180,22){$\rho_{\phi^1}$}
\put(312,22){$\rho_{\cal Q}$}
\end{picture}}
\end{equation}
\vspace{0.5ex}

\noindent
v) D'apr\`es \cite{Va}, on a un isomorphisme
$ l^\ast\colon H^{q}({\cal F},{\cal F}_{0};\nu^{\ast}{\cal F}) \rightarrow
H^{q}(M,M_{0};\phi^{1}) $.
En restriction aux fibres $F_u$, cet isomorphisme induit des morphismes
$l_u^\ast$ qui
rendent commutatif le diagramme suivant:

\begin{equation}                                \label{isophi}
\raisebox{-6ex}{
\begin{picture}(240,55)(0,0)
\put(105,46){\vector(1,0){65}}
\put(123,10){\vector(1,0){40}}
\put(18,43){$H^{2}({\cal F},{\cal F}_0;\nu^\ast{\cal F})$}
\put(60,37){\vector(0,-1){15}}
\put(0,7){$\prod_{u \in M_0}
H^{2}(F_u;\nu^\ast{\cal F} \! \mid_{\textstyle F_u})$}
\put(177,43){$H^{2}(M,M_0;\phi^{1})$}
\put(215,37){\vector(0,-1){15}}
\put(171,7){$\prod_{u \in M_0} H^{2}(F_u;\phi^{1} \! \mid_{\textstyle F_u})$}
\put(138,52){$l^\ast$}
\put(133,15){$\{l^\ast_u\}$}
\put(42,27){$\rho_{\cal F}$}
\put(223,27){$\rho_{\phi^1}$}
\end{picture}}
\end{equation}
\vspace{0.5ex}

\noindent
On remarque que le groupe
\( H^{2}(F_u;\nu^\ast{\cal F} \!\! \mid_{\textstyle F_u}) \)
est isomorphe au groupe de cohomologie de De Rham
\( H^{2}(F_u;\nu^{\ast}_{u}{\cal F}_0) \cong
H^{2}(F_u;{\Bbb R}^{m}) \).
\vspace{2ex}

\noindent
vi) Le faisceau restreint
\( {\cal P} \!\! \mid_{\textstyle F_u} \)
est le faisceau constant dont la fibre est la fibre
\( ({\cal P}_{0})_u \cong {\Bbb Z}^{r} \)
de ${\cal P}_0$ en $u$. Donc le groupe
\( H^{2}(F_u;{\cal P} \!\! \mid_{\textstyle F_u}) \)
est isomorphe au groupe de cohomologie singuli\`{e}re
\( H^{2}(F_u;({\cal P}_0)_u) \cong
H^{2}(F_u;{\Bbb Z}^{r}) \).
On a la factorisation suivante:

\begin{equation}
\raisebox{-4ex}{
\begin{picture}(190,55)(0,0)
\put(82,37){\vector(1,-1){17}}
\put(82,44){\vector(1,0){50}}
\put(112,21){\vector(1,1){15}}
\put(2,41){$ H^{2}(F_u;{\cal P} \!\! \mid_{\textstyle F_u})$}
\label{fact} \put(140,41){$H^{2}(F_u;\phi^{1} \! \mid_{\textstyle F_u})$}
\put(65,5){$H^{2}(F_u;\nu^\ast{\cal F} \!\! \mid_{\textstyle F_u})$}
\put(73,23){$h^\ast_u$}
\put(102,49){$j^\ast_u$}
\put(130,23){$l^\ast_u$}
\end{picture}}
\end{equation}
\vspace{1ex}

 L'hypoth\`ese d'int\'egrabilit\'e cohomologique est une reformulation du
crit\`ere d'int\'egrabilit\'e de P. Dazord (\cite{D2}). Le th\'eor\`eme suivant
va montrer que ce crit\`ere est effectivement v\'erifi\'e pour les
vari\'et\'es de
Poisson consid\'er\'ees:

\begin{th}             \label{thintco}
Une vari\'{e}t\'{e} de Poisson $(M_{0},\Lambda_{0})$ qui v\'erifie les
conditions
$(H_1)$ et ($H_2$) est cohomologiquement int\'{e}grable,i.e. il existe
\( \nu \in H^{2}(M,M_{0};{\cal P}) \)
telle que
\( j^{\ast}\nu = d_{1}[\Lambda] \).
\end{th}

\noindent
{\bf D\'{e}monstration}
Soient $\Omega$ un repr\'esentant de la classe $d_1[\Lambda]$ et
$\Omega_u$ sa restriction \`a la fibre $F_u$. On va montrer que la construction
d'un ant\'ec\'edent $\nu$ de $d_1[\Lambda]$ se ram\`ene \`a la
construction de classes enti\`eres
\( \nu_{u} \in H^{2}(F_u;{\cal P}\! \mid_{\textstyle F_u}) \)
qui "int\`egrent" les classes r\'eelles
\( [\Omega_u] \in H^{2}(F_u;\nu^{\ast}{\cal F} \! \mid_{\textstyle F_u}) \).
D'apr\`es l'exactitude de la suite~(\ref{coQ}),
la classe $d_1[\Lambda]$ se remonte en une classe $\nu$ si et seulement si
la classe projet\'ee $k^{\ast}(d_1[\Lambda])$ est nulle. Si l'on consid\`ere le
diagramme~(\ref{diares}), on obtient:
$$
\rho_{\cal Q}(k^{\ast}(d_1[\Lambda])) =
\{k^{\ast}_u\}(\rho_{\phi^1}(d_1[\Lambda]))
$$
Or, d'apr\`es le th\'{e}or\`{e}me~\ref{thcoh}
(que l'on prouvera dans \S6), le morphisme de restriction
$$
\rho_{\cal Q}\colon H^{2}(M,M_0;{\cal Q}) \longrightarrow
\prod_{u \in M_0}  H^{2}(F_u;{\cal Q} \! \mid_{\textstyle F_u})
$$
est injectif car les fibres $F_u$ sont simplement connexes. Il s'ensuit que
$k^{\ast}(d_1[\Lambda])=0$ si et seulement si la famille de restrictions
$$
\rho_{\phi^1}(d_1[\Lambda]) = \{ l^{\ast}_u[\Omega_u]\} \in
\prod_{u \in M_0} H^{2}(F_u;\phi^{1} \! \mid_{\textstyle F_u})
$$
(voir le diagramme~(\ref{isophi})) se projette par $\{k^\ast_u\}$ sur la classe
nulle. Pour cela, il faut et il suffit que les classes r\'eelles $[\Omega_u]$
poss\`edent des ant\'ec\'edents entiers $\nu_u$ car
$$
j^\ast_u (\nu_u) =
l^\ast_u ( h^\ast_u (\nu_u) ) =
l^{\ast}_u[\Omega_u]
$$
d'apr\`es la factorisation~(\ref{fact}).

Enfin, les fibres $F_u$ \'{e}tant simplement connexes,
l'int\'{e}gration de $\Omega$ sur les sph\`{e}res tangentes aux fibres
d\'efinit des classes
$$
\nu_{u} \in Hom(\pi_2(F_u),{\cal P} \!\! \mid_{\textstyle F_u}) =
H^{2}(F_u;{\cal P} \!\! \mid_{\textstyle F_u})
$$
telles que
\( h^\ast_u(\nu_u) = [\Omega_u] \);
d'o\`{u} le th\'{e}or\`{e}me. $\Box$
\vspace{1ex}

\subsection{Int\'egration diff\'erentiable: construction du fibr\'{e} principal
\`a groupo\"{\i}de structural}  \label{fibre}

\hspace{1em}
La notion de {\bf fibr\'e principal \`a groupo\"{\i}de structural} est
d\^ue \`a A. Haefliger (\cite{H}). Soient
\( {\cal H} \sgroup{\beta}{\alpha} {\cal H}_0 \)
un groupo\"{\i}de de Lie,
$\Gamma$ une vari\'{e}t\'{e} munie d'une application diff\'{e}rentiable
\( \rho\colon \Gamma \rightarrow {\cal H}_0 \)
et
$$
\Gamma \times_{{\cal H}_0} {\cal H} =
\{ (\gamma,h) \in \Gamma \times {\cal H} \colon \alpha(\gamma)=\rho(h) \}
$$
le produit fibr\'{e} correspondant.
Une {\bf action \`a droite de $\cal H$ sur $\Gamma$} est
une application diff\'erentiable
\( (\gamma,h)  \in  \Gamma \times_{{\cal H}_0} {\cal H}
\mapsto  \gamma.h  \in \Gamma \)
v\'erifiant des propri\'et\'es qui g\'en\'eralisent
celles de l'action d'un groupe (voir \cite{M}).
On dira que l'action est {\bf propre} si l'application
\( (\gamma,h)  \in  \Gamma \times_{{\cal H}_0} {\cal H}
\mapsto  (\gamma,\gamma.h)  \in  \Gamma \times \Gamma  \)
est propre.
De fa\c{c}on concr\`ete, une action libre et propre
d\'{e}finit une structure de $\cal H$-fibr\'{e} principal
\( \pi\colon \Gamma \rightarrow M = \Gamma/{\cal H} \).

On supposera d\'esormais que
$\cal H$ est le groupo\"{\i}de structural $\cal G$ de $(M,\Lambda)$
et l'on rappelle que celui-ci est un fibr\'e en groupes ab\'eliens
$\cal F$-feuillet\'e.

Les notions d'{\bf atlas fibr\'{e}} et {\bf
cocycle} s'\'{e}tendent de fa\c{c}on naturelle au cas
des fibr\'{e}s principaux \`{a} groupo\"{\i}de structural.
Ceux-ci sont classifi\'{e}s par la classe
\( \tau \in H^{1}(M;
\raisebox{-1.1ex}{$\stackrel{\textstyle {\cal G}}{\sim}$}) \)
d'un cocycle (voir \cite{H}).
Un $\cal G$-fibr\'{e} principal au-dessus de $M$ est trivial en
restriction \`{a} $M_{0}$ si et seulement si
le cocycle induit est cohomologue \`{a} z\'{e}ro.
Puisque la suite~(\ref{cohrel}) est exacte,
ce sera le cas si la classe $\tau$ est relative.
\vspace{1ex}

\begin{prop}       \label{fp}
L'int\'egration cohomologique
\( \nu \in H^{2}(M,M_{0};{\cal P}) \)
de $(M_0,\Lambda_0)$ est repr\'esent\'ee par un cocycle qui
d\'efinit un $\cal G$-fibr\'e principal
\( \pi\colon \Gamma \rightarrow M \).
De plus, ce fibr\'e est trivial en restriction \`a $M_0$.
\end{prop}

\noindent
{\bf D\'emonstration}
La classe
\( \tau = \delta^{-1} \nu \in H^{1}(M,M_{0};
\raisebox{-1.1ex}{$\stackrel{\textstyle {\cal G}}{\sim}$}) \)
est repr\'{e}sent\'{e}e par un co\-cycle
\( ( \{ U_i \} , \{ \tau_{ij} \} ) \)
sur $M$ \`{a} valeurs dans $\cal G$.
Soit $\Gamma$ le quotient de la r\'{e}union disjointe
\( \coprod {\cal G} \!\! \mid_{\textstyle U_i} \)
par la relation d'\'{e}quivalence qui identifie
$g$ dans
\( {\cal G} \!\! \mid_{\textstyle U_j} \)
avec $g + \tau_{ij}(p(g))$ dans
\( {\cal G} \!\! \mid_{\textstyle U_i} \).
La projection
\( p\colon \coprod {\cal G} \!\! \mid_{\textstyle U_i} \rightarrow M \)
passe au quotient en une submersion surjective
\( \pi\colon \Gamma \rightarrow M \).
La projection
\( \psi\colon \coprod {\cal G} \!\! \mid_{\textstyle U_i} \rightarrow \Gamma \)
se restreint en un diff\'{e}omorphisme $\cal G$-\'{e}quivariant
$\psi_i$ qui rend commutatif le diagramme suivant:

\begin{center}
\begin{picture}(100,43)(0,0)
\put(32,28){\vector(1,-1){17}}
\put(32,33){\vector(1,0){50}}
\put(82,28){\vector(-1,-1){17}}
\put(0,33){${\cal G} \!\! \mid_{\textstyle U_i}$}
\put(90,33){$\pi^{-1}(U_i)$}
\put(53,0){$U_i$}
\put(29,14){$p$}
\put(79,14){$\pi$}
\put(54,37){$\psi_i$}
\end{picture}
\end{center}
\vspace{0.2ex}

\noindent
Donc son inverse $\varphi_i$ est une carte de trivialit\'{e} locale.
L'atlas fibr\'{e} $\{ (U_i,\varphi_i) \}$
d\'{e}finit une structure de $\cal G$-fibr\'{e} principal sur $\Gamma$.
Puisque la classe $\tau$ est relative, ce fibr\'e est trivial en restriction
\`a $M_0$
et l'on a le diagramme suivant:

\begin{equation}                     \label{find}
\raisebox{-4ex}{
\begin{picture}(150,53)(0,0)
\put(46,42){\vector(1,0){60}}
\put(46,6){\vector(1,0){60}}
\put(26,39){${\cal G}_0$}
\put(117,39){$\Gamma$}
\put(74,48){$i$}
\put(46,13){\vector(3,1){65}}
\put(26,31){\vector(0,-1){15}}
\put(34,17){\vector(0,1){15}}
\put(122,31){\vector(0,-1){15}}
\put(22,3){$M_0$}
\put(115,3){$M$}
\put(74,10){$\varepsilon_{0}$}
\put(10,21){$p_0$}
\put(38,21){$s_0$}
\put(130,21){$\pi$}
\put(74,28){$\varepsilon$}
\end{picture}}
\end{equation}
\vspace{0.3ex}

\noindent
o\`u la section nulle $s_0$ de $p_0$ d\'efinit une section $\varepsilon$
de $\pi$ au-dessus de $M_0$. $\Box$
\vspace{1ex}

On remarquera que pour l'essentiel cette construction correspond \`a la
cons\-truction d'une
{\bf r\'ealisation isotrope de Libermann} par P. Dazord dans \cite{D2}.

\subsection{Int\'egration diff\'erentiable: construction d'une conne\-xion
$\cal F$-partielle}
\label{con}

\hspace{1em}
La projection
\(q\colon \nu^{\ast}{\cal F} \rightarrow {\cal G} \)
est l'application exponentielle
de l'alg\'{e}bro\"{\i}de de Lie $\nu^{\ast}{\cal F}$ sur $\cal G$.
A toute section $\mu$ de $\nu^{\ast}{\cal F}$, on associe le {\bf champ
fondamental}
$\mu^\ast$ d\'efini par
$$
\mu^{\ast}_{\gamma} = \frac{d}{dt} (\gamma.q(t\mu)) \!\! \mid_{t=0}
\;\;\;,\;\;\; \gamma \in \Gamma
$$

Soit $\cal S$ le feuilletage image r\'{e}ciproque de $\cal F$ par $\pi$.
Une {\bf 1-forme de connexion $\cal F$-partielle} est une
1-forme feuillet\'{e}e
\( \Theta \in \Omega^{1}({\cal S};\nu^{\ast}{\cal S}) \)
\`{a} valeurs dans l'alg\'{e}bro\"{\i}de de Lie
\( \nu^{\ast}{\cal S} \)
qui v\'{e}rifie les deux conditions
suivantes:
\vspace{1ex}

\noindent
1)
\( \Theta(\mu^{\ast}) = \pi^{\ast}\mu \:\: \)
pour toute section $\mu$ de l'alg\'{e}bro\"{\i}de de Lie
\( \nu^{\ast}{\cal F} \).
\vspace{1ex}

\noindent
2) $\Theta$ est {\bf $\cal G$-invariante},
i.e. pour toute section feuillet\'{e}e $\hat{\mu}$
de $\cal G$, on a:
$$
(R_{\textstyle \hat{\mu}})^{\ast} \Theta = \Theta
$$
o\`u $R_{\textstyle \hat{\mu}}$ est la {\bf translation \`{a} droite}
correspondante pour l'action de $\cal G$ sur $\Gamma$.

\begin{exem}
{\em
L'action
\(  (\mu_1,\mu_2) \in
\nu^{\ast}{\cal F} \times_M \nu^{\ast}{\cal F}
\mapsto
\mu_1 + \mu_2 \in \nu^{\ast}{\cal F} \)
est libre et propre. Le fibr\'{e} conormal
\( \nu^{\ast}{\cal F} \)
est donc muni d'une structure naturelle de fibr\'{e} principal
de groupo\"{\i}de structural
\( \nu^{\ast}{\cal F} \).
La connexion partielle de Bott s'interpr\`{e}te
(par l'interm\'{e}diaire de la forme de
Maurer-Cartan $\Phi$, voir \S\ref{groupat}) comme une connexion
partielle {\bf plate} (voir \S\ref{courb}).}
\end{exem}

La donn\'ee de la connexion partielle plate $\Phi$ sur $\cal G$ est
essentielle pour l'existence de connexions partielles sur $\Gamma$:
elle jouera le m\^eme r\^ole que la connexion cano\-nique des trivialisations
locales d'un fibr\'e principal classique.

\begin{prop}
Le $\cal G$-fibr\'{e} principal
\( \pi\colon \Gamma \rightarrow M \)
poss\`{e}de une 1-forme de connexion $\cal F$-partielle relative
\( \Theta \in
\Omega^{1}({\cal S},{\cal F}_0;\nu^{\ast}{\cal S}) \).
\end{prop}

\noindent
{\bf D\'{e}monstration}
Soit $\{(U_i,\varphi_i)\}$ un atlas fibr\'{e} dont le recouvrement sous-jacent
est localement fini. Soit $\{\rho_i\}$ une partition de l'unit\'{e}
subordonn\'{e}e.
La 1-forme feuillet\'{e}e
\( \Theta_i = \varphi_{i}^{\ast} \Phi \in
\Omega^{1}({\cal S} \!\! \mid_{\textstyle U_i};
\nu^{\ast}{\cal S} \!\! \mid_{\textstyle U_i}) \)
est une 1-forme de connexion partielle sur l'ouvert $\pi^{-1}(U_i)$. Alors on
v\'erifie ais\'ement que
$$
\Theta = \sum \: \pi^{\ast}(\rho_i) \: \: \Theta_{i} \: \in \:
\Omega^{1}({\cal S};\nu^{\ast}{\cal S})
$$
est une 1-forme de connexion partielle sur $\Gamma$.

Il reste \`a montrer que la 1-forme feuillet\'ee $\Theta$ est relative,
c'est-\`a-dire que $\varepsilon^{\ast}\Theta=0$ o\`u
$\varepsilon = i {\scriptstyle \circ} s_0$, voir le diagramme~(\ref{find}).
Or la 1-forme de connexion partielle induite
$i^\ast \Theta$ est la 1-forme de Maurer-Cartan $\Phi_0$ de ${\cal G}_0$ et
celle-ci est nulle en restriction \`{a} $M_0$
d'apr\`es~(\ref{prouni}). Il s'ensuit que
$\varepsilon^\ast \Theta =
s_0^\ast ( i^\ast \Theta ) = s_0^\ast \Phi_0 = 0$. $\Box$
\vspace{3ex}

En restriction \`{a} chaque ouvert $U_i$,
la diff\'{e}rence $\Theta - \Theta_i$ des connexions partielles
$\Theta$ et
\( \Theta_i = \varphi_{i}^{\ast} \Phi \)
se projette en une 1-forme feuillet\'{e}e relative $\theta_i$.
Le r\'{e}sultat suivant caract\'{e}rise
le recollement de connexions partielles locales:
\vspace{1ex}

\begin{prop}      \label{reccon}
Soit $\tau_{ij}$ le cocycle qui d\'efinit le $\cal G$-fibr\'e principal
$\Gamma$. Les 1-formes
feuillet\'{e}es relatives $\theta_i$ v\'{e}rifient
\( \theta_j - \theta_i = d_{\cal F} \tau_{ij} = \tau_{ij}^{\ast} \Phi \).
R\'eciproquement, toute famille de 1-formes feuillet\'{e}es relatives
$\theta_i$
v\'erifiant cette propri\'{e}t\'{e} d\'{e}termine une 1-forme de connexion
partielle
\( \Theta \in
\Omega^{1}({\cal S},{\cal F}_0;\nu^{\ast}{\cal S}) \)
\end{prop}

\noindent
{\bf D\'emonstration}
Pour toute intersection $U_i \cap U_j$, la diff\'{e}rence $\theta_j - \theta_i$
se rel\`eve en
une 1-forme feuillet\'ee sur $\pi ^{-1}(U_i \cap U_j)$ donn\'ee par
$$
\begin{array}{lcccl}
\pi^\ast ( \theta_j - \theta_i ) & = & \Theta_i - \Theta_j & = &
\varphi_i^{\ast} \Phi - \varphi_j^{\ast} \Phi \\
 & & & = & \varphi_j^\ast( (\varphi_i {\scriptstyle \circ} \varphi_j)^\ast
\Phi - \Phi ) \\
 & & & = & \varphi_j^\ast( (L_{\textstyle \tau_{ij}})^\ast \Phi - \Phi ) \\
 & & & = & \varphi_j^\ast( p^\ast (d_{\cal F} \tau_{ij}))   \\
 & & & = & \pi^\ast(d_{\cal F} \tau_{ij})
\end{array}
$$
d'apr\`es l'identit\'e~(\ref{Lmu}). Donc
$\theta_j - \theta_i = d_{\cal F} \tau_{ij} = \tau_{ij}^{\ast} \Phi$
d'apr\`es la propri\'et\'e universelle~(\ref{prouni}) de $\Phi$.
R\'eciproquement, si les 1-formes feuillet\'ees $\theta_i$ v\'erifient cette
condition,
les 1-formes de connexion  partielle
\( \Theta_i + \pi^{\ast}\theta_i \)
sur $\pi ^{-1}(U_i)$ se recollent en une 1-forme de connexion
partielle $\Theta$ sur $\Gamma$. $\Box$
\vspace{1ex}

\subsection{Int\'egration diff\'erentiable: construction de la courbure et
de la classe de Chern}
\label{courb}

\hspace{1em}
On d\'{e}finit la {\bf 2-forme de courbure} par
\( d_{\cal F}\Theta \in \Omega^{2}({\cal S},{\cal F}_0;\nu^{\ast}{\cal S}) \).
Comme dans le cas classique, la courbure se projette en une
2-forme feuillet\'{e}e ferm\'{e}e
\( \Omega \in
\Omega^{2}
({\cal F},{\cal F}_0;
\nu^{\ast}{\cal F}) \).

\begin{th}        \label{Chern1}
Le morphisme
\( j^{\ast}\colon H^{2}(M,M_0;{\cal P}) \rightarrow
H^{2}({\cal F},{\cal F}_0;\nu^{\ast}{\cal F}) \)
envoie la classe $\nu$ d'un cocycle sur la classe
$[\Omega]$ de la courbure d'une connexion partielle.
\end{th}

\noindent
{\bf D\'{e}monstration}
Soit $(\{U_i\},\{\tau_{ij}\})$ un cocycle qui repr\'{e}sente la classe
$\tau = \delta^{-1}\nu$ tel que le recouvrement sous-jacent soit localement
fini.
Soit $\{\rho_i\}$ une partition de l'unit\'{e} subordonn\'{e}e. On se
propose de construire un repr\'esentant de la classe $j^\ast \nu$ en
proc\'edant comme dans le cas de $S^1$-fibr\'es principaux (cf. \cite{BT}).
Pour cela, on  d\'efinit des $1$-formes  feuillet\'{e}es relatives
\( \theta_i = - \sum \: \rho_j \: \: d_{\cal F} \tau_{ij} \)
qui v\'{e}rifient
\( \theta_j - \theta_i = d_{\cal F} \tau_{ij} \). Puisque ces diff\'erences
sont
$d_{\cal F}$-ferm\'{e}es, les 2-formes feuillet\'{e}es relatives locales
$d_{\cal F}\theta_i$ se
recollent en une 2-forme feuillet\'{e}e relative globale $\Omega$
qui repr\'{e}sente la classe $j^{\ast}\nu$ (cf. \cite{BT}). Or la courbure
de la connexion
partielle $\Theta$ construite dans la proposition~\ref{reccon}
est donn\'{e}e par
\( d_{\textstyle {\cal S}} \Theta = d_{\textstyle {\cal S}} ( \Theta_i +
\pi^{\ast}\theta_i ) =  \pi^{\ast}d_{\textstyle {\cal F}}\theta_i =
\pi^{\ast}\Omega \).
Bref,
\( j^{\ast}\nu = [\Omega] \)
est la classe de la courbure $\Omega$ de la connexion partielle $\Theta$.
$\Box$
\vspace{2ex}

La classe
\( [\Omega] \in H^{2}({\cal F},{\cal F}_0;\nu^{\ast}{\cal F}) \)
sera appel\'{e}e la {\bf classe de Chern r\'{e}elle} du $\cal G$-fibr\'{e}
principal $\pi\colon \Gamma \rightarrow M$.
Le th\'eor\`eme~\ref{Chern1} justifie la terminologie suivante:
on dira que
\( \nu = \delta \tau \in H^{2}(M,M_{0};{\cal P}) \)
est la {\bf classe de Chern enti\`{e}re}
(cf. \cite{D1}, \cite{D2}).
\vspace{1ex}

\begin{cor}
La classe d\'eriv\'ee $d_{1}[\Lambda]$ est la classe de Chern r\'{e}elle de
$\Gamma$. $\Box$
\end{cor}

\subsection{Int\'egration symplectique: construction de la structure
symplectique}

\hspace{1em}
Soit
\( \sigma \in \Omega^{2}({\cal F},{\cal F}_{0}) \)
la forme feuillet\'{e}e symplectique de la structure de Poisson relev\'{e}e
$\Lambda$
sur $M$. Dans cette \'{e}tape, on cherche \`{a}
construire un repr\'{e}sentant symplectique $\eta$ de la forme
feuillet\'{e}e relev\'{e}e
\( \pi^{\ast}\sigma \in \Omega^{2}({\cal S},{\cal F}_{0}) \).
\vspace{2ex}

On fixe tout d'abord un couple de d\'{e}compositions adapt\'{e}es de
$TM_0$ et $TM$. Si les groupo\"{\i}des structuraux ${\cal G}_0$ et $\cal G$
sont s\'epar\'es,
il existe des d\'ecompositions de $T{\cal G}_0$ et $T{\cal G}$
adapt\'ees \`a celles de $TM_0$ et $TM$. Si ce n'est pas le cas, il faut
utiliser l'hypoth\`ese ($H_3)$ pour obtenir de
telles d\'ecompositions. Celles-ci d\'efinissent des d\'ecompositions
adapt\'ees de $T\Gamma$ en restriction aux ouverts de M qui trivialisent $\pi$.
Par un argument de partitions de l'unit\'e, on obtient une
d\'{e}composition de $T\Gamma$ adapt\'{e}e au couple de d\'{e}part. La
projection $\pi$ induit alors un morphisme de complexes qui rend commutatif le
diagramme suivant:

\begin{center}
\begin{picture}(220,50)(0,0)
\put(85,39){\vector(1,0){40}}
\put(85,3){\vector(1,0){40}}
\put(40,28){\vector(0,-1){15}}
\put(168,28){\vector(0,-1){15}}
\put(0,36){$\Omega^{q}({\cal F},{\cal F}_{0};\nu^{\ast}{\cal F})$}
\put(133,36){$\Omega^{q}({\cal S},{\cal F}_{0};\nu^{\ast}{\cal F})$}
\put(102,43){$\pi^{\ast}$}
\put(7,0){$\Omega^{1,q}(M,M_{0})$}
\put(140,0){$\Omega^{1,q}(\Gamma,M_{0})$}
\put(102,7){$\pi^{\ast}$}
\put(24,20){$\cong$}
\put(176,20){$\cong$}
\end{picture}
\end{center}
\vspace{0.5ex}

Soit $\omega$ le repr\'{e}sentant pur de $\sigma$.
La 2-forme feuillet\'{e}e
\( \Omega \in \Omega^{2}({\cal S},{\cal F}_0;\nu^{\ast}{\cal S}) \)
d\'etermin\'{e}e par la 3-forme pure $d_{1,0}\omega$ est la courbure d'une
connexion partielle  $\Theta$ sur $\Gamma$. Le repr\'{e}sentant pur
\( \theta \in \Omega^{1,1}(\Gamma,M_{0}) \)
de $\Theta$ v\'{e}rifie donc
$ d_{0,1} \theta = \pi^{\ast} d_{1,0}\omega $.
\vspace{1ex}

La 2-forme relative
\( \xi = \pi^{\ast}\omega - \theta \)
repr\'{e}sente $\pi^{\ast}\sigma$. Les composantes pures de $d\xi$ de
type $(0,3)$ et $(1,2)$ sont nulles.
La composante pure de type $(2,1)$ est $d_{0,1}$-ferm\'{e}e et donc
repr\'{e}sente une classe dans le terme
\( E_{1}^{2,1}({\cal S},{\cal F}_0) \)
de la suite spectrale
de Leray-Serre de la paire $(\Gamma,M_{0})$.
\vspace{1ex}

\begin{lem}     \label{alpha}
Les fibres de $\alpha\colon \Gamma \stackrel{\pi}{\longrightarrow} M
\stackrel{\alpha_0}{\longrightarrow} M_0$ sont simplement
connexes.  \end{lem}

\noindent
{\bf D\'{e}monstration} Pour tout point $u \in M_{0}$, le $\cal G$-fibr\'{e}
principal
\( \pi\colon \Gamma \rightarrow M \)
induit un fibr\'{e} principal
\( \pi\colon \alpha^{-1}(u) \rightarrow \alpha_{0}^{-1}(u) \)
dont le groupe structural ${\cal G}_{u}$ est la fibre de $\cal G$
en $u$. La connexion $\Theta$ induit une connexion usuelle sur ce fibr\'e
principal
dont la courbure est la restriction de
$\Omega$ \`a $\alpha^{-1}_0(u)$. On consid\`ere la suite exacte d'homotopie
$$
\dots \rightarrow
\pi_2(\alpha^{-1}_0(u))
\stackrel{\partial_{\ast}}{\longrightarrow}
\pi_{1}({\cal G}_{u}) = {\cal P}_{u} \longrightarrow
\pi_{1}(\alpha^{-1}(u)) \longrightarrow
\pi_{1}(\alpha^{-1}_0(u)) = 0
$$
Le bord $\partial_{\ast}$
envoie la classe d'homotopie d'une sph\`{e}re
$s$ tangente \`a $\alpha^{-1}_0(u)$ sur la p\'{e}riode
\( \int_{s} \Omega \)
d'apr\`es \cite{Ko}. Il est donc surjectif;
d'o\`u le lemme. $\Box$
\vspace{2ex}

Le lemme~\ref{alpha} implique que le terme
\( E_{1}^{2,1}({\cal S},{\cal F}_0)=0$ d'apr\`es
\cite{DH}. Il existe donc une $d_{0,1}$-primitive
\( \zeta \in \Omega^{2,0}(\Gamma,M_{0}) \)
de la partie pure de type $(2,1)$ de $d\xi$.
\vspace{1ex}

\begin{prop}          \label{symp}
La 2-forme relative
$$
\eta = \pi^{\ast}\omega - \theta - \zeta
$$
est un repr\'{e}sentant symplectique de $\pi^{\ast}\sigma$ et donc
\( \pi\colon (\Gamma,\eta) \rightarrow (M,\Lambda) \)
est un morphisme de Poisson \`{a} fibres isotropes.
\end{prop}

\noindent
{\bf D\'{e}monstration.}
La 3-forme relative
\( d\eta = d_{2,-1}\theta - d_{1,0}\zeta \)
est pure de type $(3,0)$ et $d_{0,1}$-ferm\'{e}e, c'est-\`{a}-dire que
$d\eta$ est une forme basique qui s'annule sur $M_{0}$. Donc
$\eta$ est un repr\'{e}sentant ferm\'e de $\pi^{\ast}\sigma$. Il reste \`a
prouver que $\eta$ est
\`{a} noyau nul. Soit $X$ un vecteur tangent \`{a} $\Gamma$ en $\gamma$ tel que

\begin{equation}
i_{\textstyle X}\eta =
i_{\textstyle X}\pi^{\ast}\omega -      \label{noyau}
i_{\textstyle X}\theta -
i_{\textstyle X}\zeta = 0
\end{equation}

\noindent
Pour tout vecteur vertical $Y$, on obtient
$ \Theta(Y)(X) = - i_{X}\theta(Y) = i_{X}\eta(Y) = 0 $.
Le vecteur $X$ est donc tangent \`{a} $\cal S$ et
l'identit\'{e}~(\ref{noyau}) se r\'{e}duit aux identit\'{e}s
\( i_{\textstyle X}\pi^{\ast}\omega = 0 \) et
\( i_{\textstyle X}\theta = 0 \).
Or, la premi\`ere identit\'e implique que $X$ est vertical
et donc $X=0$ d'apr\`{e}s
la deuxi\`eme identit\'e. $\Box$
\vspace{1ex}

\subsection{Int\'egration symplectique: construction de la structure de
groupo\"{\i}de symplectique}

\hspace{1em}
On se propose de finir la d\'{e}monstration du th\'{e}or\`{e}me~\ref{th2}.
D'apr\`{e}s la
proposition~\ref{symp}, les projections
\( \alpha\colon (\Gamma,\eta) \rightarrow (M_{0},\Lambda_{0}) \)
et
\( \beta\colon (\Gamma,\eta) \rightarrow (M_{0},-\Lambda_{0}) \)
forment une
{\bf paire duale stricte} (\cite{W1}). Il ne reste donc plus qu'a v\'erifier la
proposition suivante:

\begin{prop}
La paire duale stricte $(\alpha,\beta)$ est compl\`{e}te.
\end{prop}

\noindent
{\bf D\'{e}monstration}
Il suffit de montrer que
\( \pi\colon (\Gamma,\eta) \rightarrow (M,\Lambda) \)
est un morphisme de Poisson complet car la projection
\( \alpha_0\colon (M,\Lambda) \rightarrow (M_0,\Lambda_0) \)
est compl\`ete. On rappelle que $\pi$ est
{\bf complet} au sens
de \cite{DM} si pour toute 1-forme $\mu$ sur $M$ \`a
support compact, le champ $Y$
d\'{e}fini par
\( i_Y \eta = - \pi^{\ast}\mu \)
est complet. Puisque $\pi$ est un morphisme de Poisson, le champ $Y$
se projette sur
le champ complet $X$ d\'efini par
\( i_X \omega = - \mu_{0,1} \)
o\`u $\mu_{0,1}$ est la partie pure de type $(0,1)$ de $\mu$. D'autre part,
on a:
$$
\Theta(Y) = i_{\textstyle Y}\theta =
- i_{\textstyle Y}\eta + i_{\textstyle Y}\pi^\ast \omega =
\pi^{\ast} \mu - \pi^\ast \mu_{0,1} = \pi^\ast \mu_{1,0}
$$
Le champ $Y$ est donc un champ invariant par l'action de $\cal G$ qui
se projette
sur le champ complet $X$. On en d\'eduit que $Y$ est un champ complet. $\Box$
\vspace{2ex}

D'apr\`{e}s la caract\'{e}risation des groupo\"{\i}des symplectiques
(\cite{CDW}),
toute paire duale compl\`{e}te
est un groupo\"{\i}de symplectique et donc
$(\Gamma,\eta)$ r\'{e}alise
l'int\'{e}gration symplectique universelle de $(M_{0},\Lambda_{0})$ ce qui
ach\`eve la  preuve du th\'eor\`eme~\ref{th2}.

\setcounter{equation}{0}
\setcounter{th}{0}

\section{Obstruction \`a l'int\'egrabilit\'e}

\hspace{1em}
Soit
\( (\Gamma,\eta) \)
l'int\'{e}gration symplectique universelle
d'une structure de Poisson
\( \Lambda_{0} = ({\cal F}_{0},\sigma_{0}) \) sur une vari\'et\'e $M_0$.
Le but de cette section est de prouver le th\'eor\`eme~\ref{th3},
c'est-\`a-dire
de montrer que les conditions de r\'egularit\'e repertori\'ees au
th\'eor\`eme~\ref{th1}
sont bien n\'ecessaires pour pouvoir int\'egrer la structure de
Poisson $\Lambda_0$.
Ce faisant, on reactualise dans notre contexte les travaux de P. Dazord
en leurs donnant leur sens
v\'eritable qui est de repr\'esenter $d_1[\Lambda]$ comme une classe de
Chern ainsi qu'elle
a \'et\'e introduite et d\'ecrite au \S\ref{courb}.
La d\'emarche consistera \`a retrouver en ordre invers\'e les diff\'erentes
\'etapes d\'etaill\'ees \`a la section 4.

On supposera pour simplifier que
$\Pi_1({\cal F}_0)$ est s\'epar\'e m\^{e}me si le th\'eor\`eme~\ref{th3} reste
valable dans le  cas g\'{e}n\'{e}ral des vari\'{e}t\'{e}s de Poisson
r\'{e}guli\`{e}res.
\vspace{1ex}

\subsection{Extensions de groupo\"{\i}des de Lie}

\hspace{1em}
Le {\bf sous-groupo\"{\i}de d'isotropie}
\( Is\Gamma =
\{\gamma \in \Gamma \colon \alpha(\gamma)=\beta(\gamma) \} \)
de $\Gamma$ est plong\'{e}, ferm\'{e} et distingu\'{e} dans $\Gamma$.
Il n'est pas en g\'{e}n\'{e}ral \`{a} fibres connexes, mais ce sera
le cas pour la composante connexe \( {\cal N}_{0} \) de $M_{0}$
dans $Is\Gamma$. L'action naturelle \`{a} droite de ${\cal N}_{0}$
sur $\Gamma$ est donc libre et propre.
D'apr\`{e}s le crit\`{e}re de Godement (\cite{S}), le quotient
\( M = \Gamma/{\cal N}_{0} \)
est muni d'une structure de vari\'{e}t\'{e} pour laquelle
\( \pi\colon \Gamma \rightarrow M \)
est un {\bf fibr\'{e} principal de groupo\"{\i}de structural}
\( {\cal N}_{0} \).

La vari\'{e}t\'{e} quotient $M$ h\'{e}rite de $\Gamma$ une structure de
groupo\"{\i}de
de Lie \`{a} fibres simplement connexes qui s'\'{e}tale sur
\( \Pi_{1}({\cal F}_{0}) \).
D'apr\`es \cite{P}, ils sont en fait isomorphes.
Bref, on a une {\bf extension de groupo\"{\i}des de Lie}

\begin{equation}
\raisebox{-4ex}{
\begin{picture}(260,50)(0,0)
\put(53,33){\vector(2,-1){35}}
\put(50,41){\vector(1,0){50}}
\put(120,41){\vector(1,0){50}}
\put(167,33){\vector(-2,-1){35}}
\put(175,33){\vector(-2,-1){35}}
\put(107,33){\vector(0,-1){20}}
\put(113,33){\vector(0,-1){20}}
\put(33,37){${\cal N}_{0}$}
\put(178,37){$M=\Pi_{1}({\cal F}_{0})$}       \label{extg1}
\put(107,37){$\Gamma$}
\put(105,0){$M_{0}$}
\put(55,18){$p_{0}$}
\put(97,24){$\alpha$}
\put(118,24){$\beta$}
\put(140,28){$\alpha_{0}$}
\put(162,16){$\beta_{0}$}
\put(72,45){$i$}
\put(142,45){$\pi$}
\end{picture}}
\end{equation}
\vspace{1ex}

On s'int\'{e}resse \`{a} l'extension
infinit\'esimale associ\'{e}e. On rappelle que:
\vspace{1ex}

\noindent
i) la projection $\alpha_{0}$ induit un isomorphisme
\( \alpha_{0\ast}\colon {\cal L}_{M} \rightarrow {\fraktur X}({\cal F}_{0}) \)
o\`u ${\cal L}_{M}$ est l'alg\`{e}bre de Lie des champs invariants
\`{a} gauche sur $M$ (voir \cite{He});
\vspace{1ex}

\noindent
ii) l'alg\`{e}bre de Lie ${\cal L}_{\Gamma}$
des champs invariants \`{a} gauche sur $\Gamma$ est l'image du
morphisme injectif d'alg\`{e}bres de Lie
$$
\alpha^{\#} = \eta^{\#} {\scriptstyle \circ} \alpha^{\ast}\colon
\Omega^{1}(M_{0})
\longrightarrow
{\fraktur X}(\Gamma)
$$
qui, \`{a} toute 1-forme $\mu$,
associe le champ $Y=\alpha^\#$ d\'{e}fini par
\( i_{Y}\eta = - \alpha^{\ast}\mu \)
(voir \cite{CDW}).
\vspace{1ex}

Alors puisque la projection $\alpha$ est un morphisme de Poisson, on obtient un
isomorphisme de suites exactes d'alg\`{e}bres de Lie

\begin{equation}
\raisebox{-4ex}{
\begin{picture}(285,44)(0,0)
\put(0,0){$0$}
\put(0,36){$0$}
\put(16,3){\vector(1,0){16}}
\put(16,39){\vector(1,0){16}}
\put(96,39){\vector(1,0){40}}
\put(101,3){\vector(1,0){35}}
\put(58,36){${\cal L}_{{\cal N}_{0}}$}
\put(66,14){\vector(0,1){15}}
\put(40,0){$\Omega^{1}(M_{0},{\cal F}_{0})$}
\put(180,39){\vector(1,0){40}}
\put(185,3){\vector(1,0){35}}
\put(155,36){${\cal L}_{\Gamma}$}
\put(235,36){${\cal L}_{M}$}
\put(195,45){$\pi_{\ast}$}
\put(160,14){\vector(0,1){15}}        \label{isoext}
\put(240,28){\vector(0,-1){15}}
\put(142,0){$\Omega^{1}(M_{0})$}
\put(226,0){${\fraktur X}({\cal F}_{0})$}
\put(195,7){$\Lambda^{\#}$}
\put(140,18){$\alpha^{\#}$}
\put(248,18){$\alpha_{0\ast}$}
\put(258,3){\vector(1,0){16}}
\put(258,39){\vector(1,0){16}}
\put(282,0){$0$}
\put(282,36){$0$}
\end{picture}}
\end{equation}
\vspace{0.5ex}

\noindent
o\`{u} $\pi_{\ast}$ est le morphisme d'alg\`{e}bres de Lie induit par
$\pi$.

Bref, l'extension d'alg\'{e}bro\"{\i}des de Lie~(\ref{exta})
est l'extension infinit\'{e}simale associ\'{e}e \`{a}
l'extension de groupo\"{\i}des de Lie~(\ref{extg1}).
En particulier, on a le r\'{e}sultat suivant:

\begin{prop}          \label{algG}
Le fibr\'{e} conormal
\( \nu^{\ast}{\cal F}_{0} \)
est l'alg\'{e}bro\"{\i}de de Lie de
${\cal N}_{0}$. $\Box$
\end{prop}

\subsection{Action structurale, isotropie et groupo\"{\i}de
structural} \label{Sacst}

\hspace{1em}
D'apr\`{e}s la proposition~\ref{algG} et le diagramme~(\ref{isoext}),
$ \alpha^{\#}\colon \Omega^{1}(M_{0},{\cal F}_{0})
\rightarrow
{\cal L}_{\Gamma} $
est l'action infinit\'{e}simale associ\'{e}e \`{a} l'action
de ${\cal N}_{0}$ sur $\Gamma$.
Puisque le morphisme de Poisson $\alpha$ est
complet (voir \cite{CDW}), les champs verticaux $\alpha^\#\mu$ sont complets.
Si l'on note
\( \varphi_{t}^{\alpha^{\#}\mu} \)
leurs flots, l'action infinit\'{e}simale s'int\`{e}gre
en une action \`{a} droite du groupo\"{\i}de $\nu^\ast{\cal F}_0$ sur $\Gamma$

\begin{equation}          \label{acst}
\begin{array}[t]{c}
\Gamma \times_{\textstyle M_{0}} \nu^{\ast}{\cal F}_{0} \\
(\gamma,\mu)
\end{array}
\begin{array}[t]{c}
\longrightarrow \\
\mapsto
\end{array}
\begin{array}[t]{c}
\Gamma \\
\varphi_{1}^{\alpha^{\#}\mu}(\gamma)
\end{array}
\end{equation}

\noindent
qui rel\`{e}ve l'action structurale de ${\cal N}_{0}$
sur $\Gamma$ (cf. \cite{D1}).
Son {\bf isotropie} ${\cal I}_{0}$ est engendr\'ee par les 1-formes relatives
$\mu$ pour lesquelles $ \varphi_{1}^{\alpha^{\#}\mu}$ est l'identit\'e.
Par  construction, ${\cal N}_0$ est le groupo\"{\i}de quotient
$ \nu^\ast{\cal F}_0 / {\cal I}_0 $.
\vspace{1ex}

On montre ais\'ement que l'action~(\ref{acst})
est propre et donc que celle de ${\cal N}_{0}$
est libre et propre. Bref, on retrouve
la structure de ${\cal N}_{0}$-fibr\'{e} principal sur $\Gamma$ comme
une cons\'{e}quence de la compl\'{e}tion de $\alpha$.

\subsection{Construction du fibr\'{e} principal \`{a} groupo\"{\i}de
structural}

\hspace{1em}
Soit $\Lambda=({\cal F},\sigma)$ la structure de Poisson relev\'{e}e de
$\Lambda_{0}$ sur le groupo\"{\i}de d'homotopie
\( M = \Pi_{1}({\cal F}_{0}) \).
D'apr\`{e}s \cite{D2}, la projection
\( \pi\colon (\Gamma,\eta) \rightarrow (M,\Lambda) \)
est une {\bf r\'{e}alisation isotrope de Libermann}, i.e. un
morphisme de Poisson complet \`{a} fibres connexes
s\'{e}par\'{e}es et isotropes. La compl\'etion de $\pi$ entra\^{\i}ne que tout
champ vertical \( Y = \pi^{\#} \mu \)
d\'{e}fini par
\( i_{Y}\eta = - \pi^{\ast}\mu \)
est complet et son flot est not\'{e}
\( \varphi_{t}^{\mu} \).
\vspace{1ex}

L'action de
\( \nu^{\ast}{\cal F} \)
sur $\Gamma$ d\'{e}finie par

\begin{equation}
\begin{array}[t]{c}
\Gamma \times_{\textstyle M} \nu^{\ast}{\cal F} \\
(\gamma,\mu)
\end{array}
\begin{array}[t]{c}
\longrightarrow \\
\mapsto
\end{array}
\begin{array}[t]{c}
\Gamma \\
\varphi_{1}^{\mu}(\gamma)
\end{array}
\end{equation}

\noindent
est le rel\`{e}vement par $\alpha_{0}$  de l'action
\`{a} droite~(\ref{acst}) de
\( \nu^{\ast}{\cal F}_{0} \)
sur $\Gamma$.
Son isotropie $\cal I$ est l'image
r\'{e}ciproque par
$\alpha_{0}$ et $\beta_{0}$ de l'isotropie ${\cal I}_{0}$ de l'action
bilat\`{e}re de
\( \nu^{\ast}{\cal F}_{0} \)
sur $\Gamma$. Les projections $\alpha_{0}$ et $\beta_{0}$ d\'{e}finissent
alors un m\^{e}me groupo\"{\i}de de Lie image r\'{e}ciproque
\( {\cal N} = \nu^{\ast}{\cal F} / {\cal I} \).
\vspace{1ex}

L'action de $\cal N$ sur $\Gamma$ est libre et propre et donc
$\Gamma$ est muni d'une structure de
$\cal N$-fibr\'{e} principal au-dessus de $M$. Celle-ci est essentiellement
identique  \`a la structure de ${\cal N}_{0}$-fibr\'{e} principal d\'ecrite au
\S\ref{Sacst}: le groupo\"{\i}de $\cal N$ est le mod\`{e}le local commun.
\vspace{1ex}

\subsection{Construction de la connexion $\cal F$-partielle canonique}

\hspace{1em}
Soit $Y = \pi^\#\mu$ le champ vertical complet d\'efini par une 1-forme
relative $\mu$ sur $M$. La courbe int\'{e}grale de $Y$ passant par
\( \gamma \in \Gamma \)
s'\'ecrit
$$
t \in {\Bbb R} \mapsto
\varphi_{t}^{\mu}(\gamma) =
\varphi_{1}^{t\mu}(\gamma)\in \Gamma
$$
Il s'ensuit que $Y$ est le champ fondamental de l'action de $\cal N$ sur
$\Gamma$ associ\'e
\`a la section $\mu$ de l'alg\'{e}bro\"{\i}de de Lie $\nu^{\ast}{\cal F}$.

Soit $\cal S$ le feuilletage image r\'{e}ciproque de $\cal F$ par $\pi$.
A toute d\'{e}composition de $TM_{0}$,
on associe la 1-forme feuillet\'ee
\( \Theta \in \Omega^{1}({\cal S};\nu^{\ast}{\cal S}) \)
donn\'{e}e par
$$
\Theta(\pi^{\#}\mu) =
\pi^{\ast}
\mu_{1,0}
$$
o\`{u} $\mu_{1,0}$ est la partie pure de type $(1,0)$
de la 1-forme $\mu$ pour la
d\'{e}composition de $TM$ adapt\'{e}e \`{a} celle de $TM_{0}$.
\vspace{1ex}

\begin{prop} La 1-forme feuillet\'ee $\Theta$ est une 1-forme de connexion
$\cal F$-par\-tielle
relative sur $\Gamma$ que l'on appelera {\bf canonique}.
\end{prop}

\noindent
{\bf D\'emonstration} Pour toute 1-forme relative $\mu$ sur $M$, la 1-forme
$\Theta$
v\'erifie la relation $\Theta(\pi^\# \mu) = \pi^\ast \mu$.
L'invariance de $\Theta$ est obtenue en remarquant que l'on a

\begin{equation}    \label{invcon}
(\varphi_{1}^{\mu})^{\ast} \Theta - \Theta = \pi^{\ast} d_{\cal F} \mu
\end{equation}

\noindent
Si l'on pose $Y = \pi^{\#}\mu$, alors la 1-forme feuillet\'{e}e
\( i_{\textstyle Y} d_{\cal S} \Theta \)
est nulle. En effet, pour toute 1-forme $\mu_{1}$ sur $M$, on a:
$$
d_{\cal S} \Theta (Y,Y_{1}) =
i_{\textstyle Y}d_{\cal S}\Theta(Y_{1}) -
i_{\textstyle Y_{1}}d_{\cal S}\Theta(Y) - \Theta([Y,Y_{1}]) =
\pi^{\ast}(i_{\textstyle X_{1}}d\mu +
\{\mu,\mu_{1}\}) = 0
\vspace{1ex}
$$
o\`{u} le champ
\( Y_{1} = \pi^{\#}\mu_{1} \)
se projette sur le champ
\( X_{1} = \Lambda^{\#}\mu_{1} \). Il s'ensuit que
$$
L_{\textstyle Y}\Theta =
i_{\textstyle Y}d_{\cal S}\Theta + d_{\cal S}i_{\textstyle Y}\Theta =
d_{\cal S}(i_{\textstyle Y}\Theta) =
d_{\cal S}(\pi^{\ast}\mu) =
\pi^{\ast}d_{\cal F}\mu
$$
et donc
$$
\frac{d}{dt}(\varphi_{t}^{\mu})^{\ast}\Theta =
(\varphi_{t}^{\mu})^{\ast} L_{\textstyle Y} \Theta =
\pi^{\ast}d_{\cal F}\mu
$$
L'identit\'e~(\ref{invcon}) s'en d\'{e}duit par int\'{e}gration.
Enfin, puisque l'involution $\iota$ est un
isomorphisme antisymplectique, on a
\( \iota^{\ast}\Theta = - \Theta \)
et donc la 1-forme de connexion partielle
est relative, c'est-\`a-dire
$\Theta \in \Omega^2({\cal S},{\cal F}_0;\nu^\ast{\cal S})$. $\Box$
\vspace{1ex}

\begin{prop}            \label{Chern2}
La classe de Chern r\'{e}elle du $\cal N$-fibr\'{e}
principal $\Gamma$ est \'{e}gale \`{a} la classe d\'{e}riv\'{e}e
$d_{1}[\Lambda]$.
\end{prop}

\noindent
{\bf D\'{e}monstration}
Pour tout couple
\( \mu_{1} , \mu_{2} \in \Omega^{0,1}(M) \),
on a:
$$
d_{\cal F}\Theta(Y_{1},Y_{2}) =
- \Theta([Y_{1},Y_{2}]) =
- \pi^{\ast} \{\mu_{1},\mu_{2}\}_{1,0}
$$
o\`{u} $Y_{i}$ est le champ horizontal $\pi^{\#}\mu_{i}$.
La courbure de $\Theta$ est donc donn\'ee par
$$
\Omega(X_{1},X_{2}) =
- \{\mu_{1},\mu_{2}\}_{1,0} =
i_{\textstyle X_{2}}i_{\textstyle X_{1}} d_{1,0} \omega
$$
o\`{u}
\( X_{i} = \Lambda^{\#}\mu_{i} \)
et $\omega$ est le repr\'{e}sentant pur de type
$(0,2)$ de $\sigma$. $\Box$

\subsection{Obstruction \`{a} l'int\'{e}grabilit\'{e}}

\hspace{1em}
D'apr\`{e}s la proposition~\ref{Chern2}, les p\'{e}riodes sph\'{e}riques
d\'{e}riv\'{e}es de $\Lambda$ sont les p\'{e}riodes sph\'{e}riques de la
courbure $\Omega$ de la connexion partielle canonique $\Theta$. Comme dans le
cas des $S^{1}$-fibr\'{e}s principaux (voir \cite{Ko}), on en d\'eduira la
proposition suivante:

\begin{prop}       \label{periso}
L'isotropie ${\cal I}_{0}$ de l'action de
\( \nu^{\ast}{\cal F}_{0} \)
sur $\Gamma$ est \'egale au groupo\"{\i}de d\'{e}riv\'{e} ${\cal P}_{0}$.
\end{prop}

\noindent
{\bf D\'{e}monstration}
Pour tout $u \in M_{0}$, le $\cal N$-fibr\'{e} principal
\( \pi\colon \Gamma \rightarrow M \)
induit un fibr\'{e} principal
\( \pi\colon
\alpha^{-1}(u) \rightarrow \alpha_{0}^{-1}(u) \)
dont le groupe structural ${\cal N}_{u}$ est la fibre de $\cal N$
en $u$.
Toute sph\`{e}re $s$ tangente \`{a} $\alpha_{0}^{-1}(u)$
se rel\`{e}ve en un disque $D$ tangent \`{a} $\alpha^{-1}(u)$.
Son bord est la courbe int\'{e}grale $C$ passant par $u$
d'un champ $\pi^{\#}\mu$ associ\'{e} \`{a} un \'{e}l\'{e}ment
$\mu$ de l'isotropie $\cal I$.
R\'{e}ciproquement, puisque
la fibre $\alpha^{-1}(u)$ est simplement connexe, une telle
courbe est toujours le bord d'un disque qui se projette par
$\pi$ sur une sph\`{e}re tangente \`{a}
$\alpha_{0}^{-1}(u)$. Alors on a
$$
\pi^{\ast} ( \int_{s} \Omega ) =
\int_{D} d_{\cal F}\Theta = \int_{C} \Theta =
\pi^{\ast}\mu
$$
Bref, l'isotropie $\cal I$ est \'egale au groupo\"{\i}de d\'{e}riv\'{e}
$\cal P$
 au-dessus de l'espace des unit\'{e}s $M_{0}$.
Il s'ensuit que
\( {\cal I}_{0} = {\cal P}_{0}  \)
d'apr\`{e}s le lemme~\ref{perhaut}. $\Box$
\vspace{2ex}

D'apr\`es la proposition~\ref{periso}, on obtient le th\'eor\`eme suivant dont
le th\'eor\`eme~\ref{th3} est une cons\'equence imm\'ediate:
\vspace{1ex}

\begin{th}
Le groupo\"{\i}de de Lie
${\cal N}_0$ est \'egal au groupo\"{\i}de structural
${\cal G}_0$ de la vari\'et\'e de Poisson $(M_0,\Lambda_0)$. $\Box$
\end{th}

\setcounter{equation}{0}
\setcounter{th}{0}

\section{Cohomologie des submersions}

\hspace{1em}
L'objet de cette section est de d\'emontrer le r\'{e}sultat de cohomologie
utilis\'{e} dans l'int\'{e}gration cohomologique
(voir le th\'{e}or\`{e}me~\ref{thintco}).
\vspace{1ex}

\subsection{Cohomologie relative des submersions}

\hspace{1em}
Soient $M$ une vari\'{e}t\'{e}, $B$ une sous-vari\'{e}t\'{e} de $M$ et
\( \pi\colon M \rightarrow B \)
une submersion surjective.
Soient ${\cal Q}_0$ un faisceau de base $B$ et $\cal Q$ le faisceau image
r\'{e}ciproque
$\pi^{\ast}{\cal Q}_0$.

Le faisceau de Leray ${\cal H}^{q}(\pi)$ est engendr\'{e} par le
pr\'{e}faisceau
qui, \`{a} tout ouvert $V$ de $B$, associe le groupe
$H^{q}(\pi^{-1}(V),V;{\cal
Q})$; la fibre de $b$ est la limite inductive
$$
{\cal H}^{q}(\pi)(b) =
\raisebox{-3ex}
{$\stackrel{\textstyle lim}{\stackrel{\textstyle \rightarrow}
{\scriptstyle b \in V}}$} \;
H^{q}(\pi^{-1}(V),V;{\cal Q})
$$
Les morphismes induits par les inclusions de $F_b$
dans $\pi^{-1}(V)$ d\'{e}finissent
par passage \`{a} la limite un morphisme
\( \rho_b\colon {\cal H}^{q}(\pi)(b) \rightarrow H^{q}(F_b,b;{\cal Q}) \)
qui n'est en g\'{e}n\'{e}ral ni injectif, ni surjectif. En fait, $\rho_b$
est un isomorphisme si $\{ \pi^{-1}(V) \}$ est un syst\`{e}me fondamental
de voisinages de
$F_b$; c'est le cas si $F_b$ est compacte (voir \cite{G}).
On consid\`ere la {\bf suite spectrale de Leray-Serre de cohomologie
relative de la paire}
$(M,B)$ (voir \cite{G})
$$
E_{2}^{p,q}(\pi) = H^{p}(B;{\cal H}^{q}(\pi)) \Longrightarrow
H^{p+q}(M,B;{\cal Q})
$$
Le terme
\( E_{2}^{0,q}(\pi) \)
s'identifie aux sections
\( \Gamma({\cal H}^{q}(\pi)) \)
de ${\cal H}^{q}(\pi)$. En degr\'{e} $q \geq 1$,
les morphismes $\rho_b$ d\'{e}finissent donc un morphisme
$\rho_2\colon E_{2}^{0,q}(\pi) \rightarrow
\prod_{b \in B} H^{q}(F_b;{\cal Q}) $.
Si les fibres $F_b$ sont compactes, alors $\rho_2$ est injectif.
Par construction, le diagramme suivant est commutatif:

\begin{equation}
\raisebox{-4.5ex}{
\begin{picture}(220,55)(0,0)
\put(0,0){$0$}
\put(8,3){\vector(1,0){15}}
\put(100,45){\vector(1,0){40}}
\put(85,3){\vector(1,0){60}}
\put(56,34){\vector(0,-1){15}}        \label{dia1}
\put(88,10){\vector(2,1){50}}
\put(178,18){\vector(0,1){15}}
\put(26,42){$H^{q}(M,B;{\cal Q})$}
\put(148,42){$\prod_{b \in B} H^{q}(F_b;{\cal Q})$}
\put(118,49){$\rho$}
\put(34,0){$E_{\infty}^{0,q}(\pi)$}
\put(160,0){$E_{2}^{0,q}(\pi)$}
\put(186,25){$\rho_2$}
\put(102,27){$\rho_\infty$}
\end{picture}}
\end{equation}
\vspace{0.5ex}

\noindent
o\`{u} le morphisme $\rho$ est induit par les inclusions des fibres $F_b$
dans $M$.
Si celles-ci sont compactes, alors les morphismes $\rho_2$ et $\rho_\infty$
sont
injectifs.
\vspace{1ex}

\begin{th}        \label{thcoh}
Si les fibres $F_b$ sont connexes et
$H^{1}(F_b;{\Bbb Z})=0$, alors
\vspace{1ex}

\noindent
i) $H^{q}(M,B;{\cal Q})=0$ pour $q=0,1$;
\vspace{1ex}

\noindent
ii)
\( H^{2}(M,B;{\cal Q}) = E_{\infty}^{0,2}(\pi) \)
et le morphisme
$$
\rho\colon H^{2}(M,B;{\cal Q}) \rightarrow \prod_{b \in B}
H^{2}(F_{b};{\cal Q})
$$
est injectif.
\end{th}

En degr\'{e} $0$, on consid\`{e}re le morphisme
\( \pi^{\ast}\colon H^{0}(B;{\cal Q}_0)
\longrightarrow
H^{0}(B;{\cal Q}) \).
Puisque les fibres de $\pi$ sont connexes et celles du faisceau
${\cal Q}_0$ sont discr\`{e}tes,
toute section de $\cal Q$ est constante en restriction \`{a} chaque fibre.
Donc le
morphisme $\pi^{\ast}$ qui est \'{e}videmment injectif est aussi
surjectif.
En degr\'{e} $1$ et $2$, on proc\'{e}dera en quatre \'{e}tapes.
L'\'{e}tude d'une
famille de  voisinages propres de $B$ dans la 1\`{e}re \'{e}tape permettra
de construire dans la 2\`{e}me \'{e}tape un recouvrement $\{W_n\}$
{\bf annihilant $H^{1}(M,B;{\cal Q})$}.
Dans la 3\`{e}me \'{e}tape, on v\'{e}rifiera l'annulation de certaines
termes de la suite spectrale de \v{C}ech de $\{W_n\}$. La preuve
du th\'{e}or\`{e}me combinera ce r\'{e}sultat avec la comparaison des suites
spectrales de Leray-Serre des restrictions $\pi_n = \pi \!\! \mid_{W_n}$.
\vspace{1ex}

\subsection{Tubes de $B$}

\hspace{1em}
On fixe une m\'{e}trique riemannienne sur $M$ dont la m\'{e}trique induite sur
chaque fibre $F_b$ est compl\`{e}te. \`{A} toute fonction continue strictement
positive
\( r\colon B \rightarrow {\Bbb R} \),
on associe la r\'{e}union $W(r)$ des boules g\'{e}od\'{e}siques ferm\'{e}es
dans
$F_b$ de centre $b$ et de rayon $r(b)$. Le voisinage $W(r)$ de $B$ sera
appel\'e le {\bf tube de $B$ de rayon $r$}.
On montre tout d'abord deux propri\'et\'es des tubes:
\vspace{1ex}

\begin{lem} \label{lem1}
Pour tout tube $W=W(r)$, la projection
\( \pi \!\! \mid_W\colon W \rightarrow B \)
est propre.
\end{lem}

\noindent
{\bf D\'{e}monstration}
Pour tout point $b \in B$, soit $U_b$ un voisinage compact de la boule
g\'{e}od\'{e}sique $W_b$ dans $F_b$. La fonction rayon $r$ \'{e}tant continue,
il existe un voisinage ouvert relativement compact $V_b$ de $b$ et
un voisinage produit
\( \overline{V_b} \times U_b \)
de $U_b$ (o\`u l'on note $\overline{V_b}$ l'adh\'erence de $V_b$) dont
l'intersection avec $W$ est un compact satur\'{e} pour  $\pi \!\! \mid_W$.
Pour tout compact $K \subset B$, on extrait un sous-recouvrement fini
\( \{ V_i \cap K \}_{i=1}^n \)
du recouvrement
\( \{ V_b \cap K \} \).
Alors
$$
(\pi \!\! \mid_W )^{-1} (K) = \bigcup_{i=1}^{n} ((\overline{V_i} \cap K)
\times U_{i}) \cap W
$$
est compact. $\Box$
\vspace{1ex}

\begin{lem}  \label{lem2}
Pour tout voisinage $U$ de $B$ tel que la projection
\( \pi \!\! \mid_U\colon U \rightarrow B \)
est propre, il existe un tube $W$ tel que $U \subset W$.
\end{lem}

\noindent
{\bf D\'{e}monstration}
Puisque $\pi \!\! \mid_U$ est propre, toute fibre $U_b$ est compacte et
donc il existe une boule g\'{e}od\'{e}sique $W_b$ dont
l'int\'{e}rieur $\stackrel{\circ}{W_b} \; \supset U_b$. Du recouvrement ouvert
\( \{ V_b \times \stackrel{\circ}{W_b} \} \),
on extrait un rafinement localement fini
\( \{ V_i \times \stackrel{\circ}{W_i} \} \).
Soit $r_i$ le rayon de la boule $W_i$. La fonction
\( s: B \rightarrow {\Bbb R} \)
d\'{e}finie par
\( s(b) = sup \; \{ \; r_i \; \colon b \in W_i \; \} \)
est une fonction en escalier. Si
\( r \!\colon \! B \rightarrow {\Bbb R} \)
est une fonction continue telle que
\( r(b) \geq s(b) \)
pour tout $b \in B$, alors $U \subset W=W(r)$. $\Box$
\vspace{2ex}

Enfin, on s'int\'{e}resse aux propri\'{e}t\'{e}s de la suite spectrale de
$\pi \!\! \mid_W$:
\vspace{1ex}

\begin{lem}  \label{lem3}
Si $W$ est un tube de $B$; alors on a
\vspace{1ex}

\noindent
i) $H^{0}(W,B;{\cal Q})=0$;
\vspace{1ex}

\noindent
ii)
\( H^{1}(W,B;{\cal Q}) = E_{2}^{0,1}(\pi \!\! \mid_W) \);
\vspace{1ex}

\noindent
iii) les morphismes
$$
\rho\colon H^{1}(W,B;{\cal Q}) \longrightarrow
\prod_{b \in B} H^{1}(W_b;{\cal Q})
$$
et
$$
\rho_\infty\colon E_{\infty}^{0,2}(\pi \!\! \mid_W)
\longrightarrow \prod_{b \in B} H^{2}(W_b;{\cal Q})
$$
sont injectifs.
\end{lem}

\noindent
{\bf D\'{e}monstration}
Les fibres $W_b$ \'etant connexes, le faisceau de Leray
\( {\cal H}^{0}(\pi \!\! \mid_ W) = 0 \).
Par cons\'{e}quent, on a:

\noindent
i) $H^{0}(W,B;{\cal Q})=0$;

\noindent
ii)
\( H^{1}(W,B;{\cal Q}) = E_{2}^{0,1}(\pi \!\! \mid_W) \).

\noindent
iii) Puisque
\( \pi \! \! \mid_W \)
est propre, le morphisme
$ \rho_2\colon E_{2}^{0,q}(\pi \!\! \mid_W)
\rightarrow \prod_{b \in B} H^{q}(W_b;{\cal Q}) $
est injectif en degr\'e $1$ et $2$. Donc les morphismes
$ \rho\colon H^{1}(W,B;{\cal Q}) \rightarrow \prod_{b \in B}
H^{1}(W_b;{\cal Q}) $
et
$ \rho_\infty\colon E_{\infty}^{0,2}(\pi \!\! \mid_W)
\longrightarrow \prod_{b \in B} H^{2}(W_b;{\cal Q}) $
sont aussi injectifs.
$\Box$
\vspace{2ex}

\subsection{Recouvrement d'annihilation}

\vspace{2ex}

\hspace{1em}
On se propose de construire un {\bf recouvrement annihilant
$H^{1}(M,B;{\cal Q})$}, i.e. un recouvrement ferm\'{e} croissant
$\{ W_n \}$ de $M$ tel que pour tout $n$,
\vspace{1ex}

\noindent
i) $W_n$ est un tube et l'application
\( \pi_n\colon W_n \rightarrow B \)
est propre \`{a} fibres connexes;
\vspace{1ex}

\noindent
ii) le morphisme de restriction
\( H^{1}(W_{n+1},B;{\cal Q}) \rightarrow H^{1}(W_n,B;{\cal Q}) \)
est nul;
\vspace{1ex}

\noindent
iii) le morphisme de faisceaux
\( {\cal H}^{1}(\pi_{n+1}) \rightarrow {\cal H}^{1}(\pi_n) \)
est nul.
\vspace{1ex}

\noindent
Pour cela, on montre tout d'abord la condition d'annihilation suivante:

\begin{lem}     \label{lemres}
Soient $W \subset \widehat{W}$ deux tubes de $B$. Si pour point tout
$b \in B$, le
morphisme de restriction
$ r_b\colon H^{1}(\widehat{W}_b;{\cal Q}) \rightarrow H^{1}(W_b;{\cal Q}) $
est nul, alors le morphisme de restriction
$$
r\colon H^{1}(\widehat{W},B;{\cal Q}) \rightarrow H^{1}(W,B;{\cal Q})
$$
et le morphisme de faisceaux
$$
r\colon {\cal H}^{1}(\pi \!\! \mid_{\widehat{W}}) \rightarrow
{\cal H}^{1}(\pi \!\! \mid_W)
$$
sont nuls.
\end{lem}

\noindent
{\bf D\'{e}monstration}
On consid\`{e}re le diagramme commutatif suivant:

\begin{center}
\begin{picture}(220,54)(0,0)
\put(97,45){\vector(1,0){45}}
\put(93,6){\vector(1,0){50}}
\put(56,34){\vector(0,-1){15}}
\put(178,34){\vector(0,-1){15}}
\put(20,42){$H^{1}(\widehat{W},B;{\cal Q})$}
\put(151,3){$\prod_{b \in B} H^{1}(W_b;{\cal Q})$}
\put(118,49){$r$}
\put(110,11){$\{r_b\}$}
\put(7,3){$\prod_{b \in B} H^{1}(\widehat{W}_b;{\cal Q})$}
\put(150,42){$H^{1}(W,B;{\cal Q})$}
\put(186,24){$\rho$}
\put(40,24){$\widehat{\rho}$}
\end{picture}
\end{center}

\noindent
D'apr\`{e}s le lemme~\ref{lem3}, le morphisme
$\rho$ est
injectif. On en d\'{e}duit que le morphisme $r$ est nul car
$\{r_b\}$ est nul. Par ailleurs,
pour tout $b \in B$, on a un diagramme

\begin{center}
\begin{picture}(220,54)(0,0)
\put(92,48){\vector(1,0){45}}
\put(90,6){\vector(1,0){55}}
\put(56,35){\vector(0,-1){17}}
\put(178,35){\vector(0,-1){17}}
\put(17,45){${\cal H}^{1}(\pi \!\! \mid_{\widehat{W}})(b)$}
\put(153,3){$H^{1}(W_b;{\cal Q})$}
\put(112,52){$r_b$}
\put(115,10){$r_b$}
\put(23,3){$H^{1}(\widehat{W}_b;{\cal Q})$}
\put(148,45){${\cal H}^{1}(\pi \!\! \mid_W)(b)$}
\put(186,26){$\rho_b$}
\put(36,26){$\widehat{\rho}_b$}
\end{picture}
\end{center}

\noindent
o\`{u} $\rho_b$ est un isomorphisme.
Par cons\'{e}quent, le morphisme de faisceux $r$
est nul fibre \`{a} fibre et donc il est nul. $\Box$
\vspace{1ex}J

\begin{prop}    \label{restube}
Si $H^{1}(F_b;{\Bbb Z})=0$, alors pour tout tube $W$ de $B$,
il existe un tube $\widehat{W} \supset W$ de $B$ tel que le morphisme de
restriction
$$
r\colon H^{1}(\widehat{W},B;{\cal Q}) \rightarrow H^{1}(W,B;{\cal Q})
$$
et le morphisme de faisceaux
$$
r\colon {\cal H}^{1}(\pi \!\! \mid_{\widehat{W}}) \rightarrow
{\cal H}^{1}(\pi \!\! \mid_W)
$$
sont nuls.
\end{prop}

\noindent
{\bf D\'{e}monstration}
Pour tout point $b \in B$, soit $U_b$ un voisinage compact de la boule
g\'{e}od\'{e}sique $W_b$ dans la fibre $F_b$. Puisque $H^{1}(F_b;{\Bbb Z})=0$,
il existe
un voisinage compact $\widehat{U}_b$ de $U_b$ tel que le morphisme de
restriction
\( r\colon H^{1}(\widehat{U}_b;{\Bbb Z}) \rightarrow H^{1}(U_b;{\Bbb Z}) \)
est nul. La continuit\'{e} du rayon entra\^{\i}ne
l'existence d'un voisinage ouvert rela\-tivement compact $V_b$ de $b$ et
d'un voisinage compact $\overline{V}_b \times \widehat{U}_b$ de
$\widehat{U}_b$ tel que
\( ( \overline{V}_b \times U_b ) \cap W$ est un compact satur\'{e} pour
$\pi \!\! \mid_W$.  Soit $\{ V_i \times \stackrel{\circ}{U_i} \}$ un rafinement
localement fini du  recouvrement ouvert
$\{ V_b \times \stackrel{\circ}{U_b} \}$
de $W$. Si l'on pose
$$
\widehat{U} = \bigcup \overline{V}_i \times \widehat{U}_i \;
\supset \;  U = \bigcup \overline{V}_i \times U_i \; \supset \; W
$$
alors la projection
\( \pi \!\! \mid_{\widehat{U}}\colon \widehat{U} \rightarrow B \)
est propre. En outre, tout morphisme de restriction
\( H^{1}(\widehat{U}_b;{\cal Q}) \rightarrow H^{1}(W_b;{\cal Q}) \)
est nul, car il se factorise \`{a} travers $H^{1}(U_b;{\cal Q})$.
D'apr\`{e}s le
lemme~\ref{lem2}, il existe un tube $\widehat{W} \supset \widehat{U}$.
\'{E}videmment les morphismes
de restriction
\( r_b\colon H^{1}(\widehat{W}_b;{\cal Q}) \rightarrow H^{1}(W_b;{\cal Q}) \)
sont encore nuls, c'est-\`{a}-dire $\widehat{W}$ v\'{e}rifie la condition
du lemme~\ref{lemres}. D'o\`{u} la proposition. $\Box$
\vspace{1ex}

\begin{th}
Si $H^{1}(F_b;{\Bbb Z})=0$, alors il existe un recouvrement
${\cal W} = \{ W_n \}$ annihilant $H^1(M,B;{\cal Q})$.
\end{th}

\noindent
{\bf D\'{e}monstration}
On d\'{e}finit $W_0=B$, puis on proc\`{e}de par r\'ecurrence.
On suppose construit un tube
$W_n$ contenant le tube $W(n)$ de rayon $n$. Soit $U$ le voisinage propre
$W_n \cup W(n+1)$ de $B$. D'apr\`{e}s le lemme~\ref{lem2}, il existe un tube
$W \supset U$. On d\'{e}finit $W_{n+1}$ comme le tube $\widehat{W}$
associ\'{e} au tube $W$ d'apr\`es la proposition~\ref{restube}. $\Box$

\subsection{Suite spectrale de \v{C}ech de $\cal W$}

\hspace{1em}
 Soit ${\cal H}^{q}$ le syst\`{e}me de coefficients sur le nerf
du recouvrement ${\cal W} = \{ W_n \}$ qui, \`{a} tout p-simplexe singulier
$n_0 < \dots < n_p$, associe le groupe
$$
H^{q}(W_{n_0 \dots n_p},B;{\cal Q}) = H^{q}(W_{n_0},B;{\cal Q})
$$
o\`{u}
\( W_{n_0 \dots n_p} = W_{n_0} \cap \dots \cap W_{n_p} \).
Soit
$ \; E_{2}^{p,q}({\cal W}) \Rightarrow H^{p+q}(M,B;{\cal Q}) \; $
la {\bf suite spectrale de \v{C}ech attach\'ee au recouvrement} $\cal W$
(\cite{G}) qui v\'{e}rifie:
$$
\begin{array}{lll}
E_{1}^{p,q}({\cal W}) = & {\cal C}^p({\cal W},{\cal H}^q) = &
\prod_{n_0 < \dots < n_p} H^{q}(W_{n_0 \dots n_p},B;{\cal Q}) \\ \\
E_2^{p,q}({\cal W}) = & H^p({\cal W},{\cal H}^q) &
\end{array}
$$

\begin{th} Soit ${\cal W}$ un recouvrement annihilant
$H^{1}(M,B;{\cal Q})$. Alors on a:
$$
\begin{array}{ll}
E_1^{p,0}({\cal W}) = 0 \;\;\;\; \mbox{pour tout $p$} \\ \\
E_2^{1,1}({\cal W}) = 0
\end{array}
$$
\end{th}

\noindent
{\bf D\'{e}monstration}
D'apr\`es le lemme~\ref{lem3}, $H^0(W_n,B;{\cal Q})=0$
pour tout $n$ et donc ${\cal H}^0=0$. Il
s'ensuit que
\( E_1^{p,0}({\cal W}) = 0 \)
pour tout $p$.
D'autre part, il faut montrer que la suite
$$
E_1^{0,1}({\cal W}) = {\cal C}^0({\cal W},{\cal H}^1)
\stackrel{d_1}{\longrightarrow}
E_1^{1,1}({\cal W}) = {\cal C}^1({\cal W},{\cal H}^1)
\stackrel{d_1}{\longrightarrow}
E_1^{2,1}({\cal W}) = {\cal C}^2({\cal W},{\cal H}^1)
$$
est exacte, o\`u $d_1$ est la diff\'erentielle de \v{C}ech $\delta$.
Une 1-cocha\^{\i}ne
\( \nu  \in C^1({\cal W},{\cal H}^1) \)
associe \`a tout couple $n<m$ une classe
$$
\nu_{nm} \in H^1(W_{nm},B;{\cal Q}) = H^1(W_n,B;{\cal Q})
$$
La 1-cocha\^{\i}ne $\nu$ est un 1-cocycle si pour tout triplet $n<m<p$, la
classe
$$
\nu_{nm} - \nu_{np} + \nu_{mp} = 0
$$
en restriction \`a $W_{nmp} = W_n$. Puisque $\cal W$ est un
recouvrement d'annihilation, la classe $\nu_{mp} \in H^1(W_m,B;{\cal Q})$ est
nulle en restriction \`a $W_n$. Bref, $\nu$ est un 1-cocycle si pour tout
triplet $n<m<p$, on a:
$$
\nu_{nm} = \nu_{np}
$$
en restriction \`a $W_n$.
Soit $\mu \in C^0({\cal W},{\cal H}^1)$
la 0-cocha\^{\i}ne qui, \`a tout $n$, associe la classe
$$
\mu_n = - \nu_{nm} \in H^1(W_n,B;{\cal Q})
$$
o\`u $n<m$. Cette classe est ind\' ependante du choix de $m$. Il ne reste
plus qu'\`a v\'erifier que $\delta \mu = \nu$. Or, on a
\( \mu_m - \mu_n = - \mu_n = \nu_{nm} \)
en restriction \`a $W_{nm}=W_n$ car la restriction de $\mu_m$ \`a
$W_n$ est nulle. $\Box$
\vspace{1ex}

\begin{cor}  \label{cor1}
$$
\begin{array}{ll}
H^1(M,B;{\cal Q}) = E_2^{0,1}({\cal W}) \\ \\
H^2(M,B;{\cal Q}) = E_{\infty}^{0,2}({\cal W}) \;\;\;\;\;\; \Box
\end{array}
$$
\end{cor}
\vspace{0.5ex}

\begin{cor}     \label{cor2}
(i)
\( \;\; H^1(M,B;{\cal Q}) = 0 \);

\noindent
(ii) le morphisme de restriction
$ \{ r_n \}\colon H^2(M,B;{\cal Q}) \rightarrow \prod_{n \in {\Bbb N}}
H^2(W_n,B;{\cal Q}) $ est injectif.
\end{cor}

\noindent
{\bf D\'emonstration}
D'apr\`es le corollaire~\ref{cor1}, en degr\'e $1$ et $2$
l'inclusion
$$
0 \to E_{\infty}^{0,q}({\cal W}) \longrightarrow E_1^{0,q}({\cal W}) =
\prod_{n \in {\Bbb N}} H^q(W_n,B;{\cal Q})
$$
co\"{\i}ncide avec le morphisme de restriction
$$
\{ r_n \}\colon H^q(M,B;{\cal Q}) \rightarrow \prod_{n \in {\Bbb N}}
H^q(W_n,B;{\cal Q})
$$
qui est donc injectif. En degr\'e $1$, ce morphisme est en plus nul,
car $\cal W$ est un recouvrement d'annihilation. $\Box$
\vspace{1ex}

\subsection{Comparaison des suites spectrales}
\vspace{1ex}

\begin{prop}  \label{propE}
\( \;\;\; H^2(M,B;{\cal Q}) = E_{\infty}^{0,2}(\pi) \)
\end{prop}

\noindent
{\bf D\'emonstration} Les inclusions $W_n \subset W_{n+1} \subset M$
induisent des morphismes

\begin{equation}  \label{ressuite}
\begin{array}{cccllc}
 & \hspace{3mm} & E_2^{p,q}(\pi) & = H^p(B;{\cal H}^p(\pi)) &
\Longrightarrow &
H^{p+q}(M,B;{\cal Q}) \\
 & \hspace{3mm} & \downarrow & & & \downarrow \\
& \hspace{3mm} & E_2^{p,q}(\pi_{n+1}) & = H^p(B;{\cal H}^p(\pi_{n+1})) &
\Longrightarrow &
H^{p+q}(W_{n+1},B;{\cal Q}) \\
& \hspace{3mm} & \downarrow & & & \downarrow \\
 & \hspace{3mm} & E_2^{p,q}(\pi_n) & = H^p(B;{\cal H}^p(\pi_n)) &
\Longrightarrow &
H^{p+q}(W_n,B;{\cal Q})
\end{array}
\end{equation}
\vspace{0.5ex}

\noindent
D'autre part, les faisceaux de Leray v\'erifient:

\noindent
1)
\( {\cal H}^0(\pi) = {\cal H}^0(\pi_{n}) = 0 \; \)
car les fibres sont connexes;

\noindent
2) ${\cal H}^1(\pi) = 0$ d'apr\`es le corollaire~\ref{cor2};

\noindent
3) le morphisme de restriction
\( {\cal H}^1(\pi_{n+1}) \rightarrow {\cal H}^1(\pi_{n}) \)
est nul car $\cal W$ est un recouvrement d'annihilation.

\noindent
On en d\'eduit que
$$
\begin{array}{l}
E_2^{p,0}(\pi) = E_2^{p,0}(\pi_n) = 0 \\
E_2^{p,1}(\pi) = 0
\end{array}
$$
et le morphisme de restriction
$$
E_2^{p,1}(\pi_{n+1}) \longrightarrow E_2^{p,1}(\pi_n)
$$
est nul pour tout $p$ et tout $n$. En degr\'e $2$, le
diagramme~(\ref{ressuite}) se r\'eduit donc au diagramme suivant:

\begin{equation}  \label{ann}
\raisebox{-9ex}{
\begin{picture}(315,106)(0,0)
\put(20,0){$0$}
\put(30,3){\vector(1,0){10}}
\put(20,40){$0$}
\put(30,43){\vector(1,0){10}}
\put(102,43){\vector(1,0){20}}
\put(102,3){\vector(1,0){20}}
\put(45,40){$E_{\infty}^{1,1}(\pi_{n+1})$}
\put(70,32){\vector(0,-1){15}}
\put(58,22){$0$}
\put(50,0){$E_{\infty}^{1,1}(\pi_n)$}
\put(211,43){\vector(1,0){20}}
\put(211,3){\vector(1,0){20}}
\put(129,40){$H^{2}(W_{n+1},B;{\cal Q})$}
\put(239,40){$E_{\infty}^{0,2}(\pi_{n+1})$}
\put(161,32){\vector(0,-1){15}}
\put(271,32){\vector(0,-1){15}}
\put(137,0){$H^{2}(W_n,B;{\cal Q})$}
\put(245,0){$E_{\infty}^{0,2}(\pi_n)$}
\put(295,3){\vector(1,0){10}}
\put(295,43){\vector(1,0){10}}
\put(312,0){$0$}
\put(312,40){$0$}
\put(20,82){$0$}
\put(30,85){\vector(1,0){10}}
\put(102,85){\vector(1,0){20}}
\put(52,82){$E_{\infty}^{1,1}(\pi)$}
\put(70,74){\vector(0,-1){15}}
\put(211,85){\vector(1,0){20}}
\put(137,82){$H^{2}(M,B;{\cal Q})$}
\put(245,82){$E_{\infty}^{0,2}(\pi)$}
\put(161,74){\vector(0,-1){15}}
\put(271,74){\vector(0,-1){15}}
\put(295,85){\vector(1,0){10}}
\put(312,82){$0$}
\put(68,98){${\scriptstyle \parallel}$}
\put(68,106){$0$}
\end{picture}}
\end{equation}
\vspace{0.5ex}

\noindent
D'o\`u la proposition. $\Box$
\vspace{1ex}

\begin{prop}     \label{prof}
Le morphisme
\( \rho\colon H^2(M,B;{\cal Q}) \rightarrow \prod_{b \in B} H^2(F_b;{\cal Q})
\)
est injectif.
\end{prop}

\noindent
{\bf D\'emonstration}
On remarque tout d'abord que le morphisme $\rho$ co\"{\i}ncide avec le
morphisme
\( \rho_\infty\colon E_{\infty}^{0,2}(\pi) \rightarrow \prod_{b \in B}
H^2(F_b;{\cal Q}) \)
d'apr\`es la proposition~\ref{propE}. On consid\`ere le diagramme commutatif
suivant:

\begin{equation}
\raisebox{-4ex}{
\begin{picture}(220,52)(0,0)
\put(0,0){$0$}
\put(8,3){\vector(1,0){15}}
\put(90,45){\vector(1,0){50}}
\put(90,3){\vector(1,0){50}}
\put(60,32){\vector(0,-1){15}}       \label{dia2}
\put(178,32){\vector(0,-1){15}}
\put(36,42){$E_{\infty}^{0,2}(\pi)$}
\put(148,42){$\prod_{b \in B} H^{2}(F_b;{\cal Q})$}
\put(108,49){$\rho_\infty$}
\put(108,10){$\rho_\infty^{n+1}$}
\put(34,0){$E_{\infty}^{0,2}(\pi_{n+1})$}
\put(148,0){$\prod_{b \in B} H^{2}((W_{n+1})_b;{\cal Q})$}
\end{picture}}
\end{equation}
\vspace{0.5ex}

\noindent
o\`u $\rho_\infty^{n+1}$ est injectif d'apr\`es le lemme~\ref{lem3}.

Soit
\( \nu \! \in \! H^2(M,B;{\cal Q}) \)
une classe telle que $\rho(\nu)=0$. La restriction
$\nu_{n+1}$ de $\nu$ \`a $W_{n+1}$ se projette sur une classe de
$E_{\infty}^{0,2}(\pi_{n+1})$ dont l'image par $\rho_\infty^{n+1}$
est nulle d'apr\`es le diagramme~(\ref{dia2}). Puisque $\rho_\infty^{n+1}$
est injectif, la projection de $\nu_{n+1}$ dans
$E_{\infty}^{0,2}(\pi_{n+1})$ est nulle. La classe $\nu_{n+1}$ se remonte
donc en une classe dans $E_{\infty}^{1,1}(\pi_{n+1})$ d'apr\`es le
diagramme~(\ref{ann}). L'annulation du morphisme de restriction
\( E_{\infty}^{1,1}(\pi_{n+1}) \rightarrow E_{\infty}^{1,1}(\pi_{n}) \)
implique que la restriction $\nu_n$ de $\nu$ \`a $W_n$ est nulle.
En r\'esum\'e, la classe $\nu$ appartient au noyau du morphisme de restriction
$\{ r_n \}$. Or, d'apr\`es le corollaire~\ref{cor2}, ce morphisme est injectif
et donc
$\nu=0$. $\Box$
\vspace{2ex}

Le corollaire~\ref{cor2} et la proposition~\ref{prof} d\'emontrent donc le
th\'eor\`eme~\ref{thcoh}. On finit par une g\'en\'eralisation de ce
th\'eor\`eme en
degr\'e arbitraire.
\vspace{1ex}

\subsection{Cohomologie des submersions}

\hspace{1em}
Soient $\pi\colon M \rightarrow B$ une submersion surjective,
${\cal Q}_0$ un faisceau de base $B$ et $\cal Q$ le faisceau image
r\'eciproque $\pi^{\ast}{\cal Q}_0$.
Pour \'etendre le th\'eor\`eme~\ref{thcoh} en
degr\'e plus grand, on remplace
la section globale par une famille de sections locales d\'efinies
sur les adh\'erences des ouverts d'un recouvrement relativement compact
et localement fini de $B$. On construit alors un recouvrement
ferm\'e croissant $\{W_n\}$ de $M$ tel que
\vspace{1ex}

\noindent
i) $W_n$ est un tube de la famille de sections locales et
\( \pi_n\colon W_n \rightarrow B \)
est une application propre \`{a} fibres connexes;
\vspace{1ex}

\noindent
ii) le diagramme suivant est commutatif:
$$
\begin{picture}(340,60)(0,0)
\put(143,48){\vector(1,0){30}}
\put(143,8){\vector(1,0){30}}
\put(155,13){$\pi_n^{\ast}$}
\put(149,54){$\pi_{n+1}^{\ast}$}
\put(25,45){$E_{2}^{1,0}(\pi_{n+1}) = H^1(B;{\cal Q}_0)$}
\put(50,25){$\parallel$}
\put(30,5){$E_{2}^{1,0}(\pi_n) = H^1(B;{\cal Q}_0)$}
\put(249,48){\vector(1,0){20}}
\put(249,8){\vector(1,0){20}}
\put(181,45){$H^{1}(W_{n+1};{\cal Q})$}
\put(277,45){$E_{\infty}^{0,1}(\pi_{n+1})$}
\put(212,37){\vector(0,-1){15}}
\put(303,37){\vector(0,-1){15}}
\put(311,27){$0$}
\put(189,5){$H^{1}(W_n;{\cal Q})$}
\put(285,5){$E_{\infty}^{0,1}(\pi_n)$}
\end{picture}
$$

De fa\c{c}on pr\'ecise, on obtient le th\'eor\`eme suivant:
\vspace{1ex}

\begin{th}         \label{subm}
Si les fibres $F_b$ sont connexes et
$H^{q}(F_b;{\Bbb Z})=0$ en degr\'e $q < r$, alors
\vspace{1ex}

\noindent
i)
\( \;\; \pi^{\ast}\colon H^q(B;{\cal Q}_0) \rightarrow
H^{q}(M;{\cal Q}) \;\; \) est isomorphisme en degr\'e $q < r$;
\vspace{1ex}

\noindent
ii) la suite
$$
0 \rightarrow
H^{r}(B;{\cal Q}_0) \stackrel{\pi^{\ast}}{\longrightarrow}
H^{r}(M;{\cal Q}) \stackrel{\rho}{\longrightarrow}
\prod_{b \in B} H^{r}(F_{b};{\cal Q})
$$
est exacte. $\Box$
\end{th}
\vspace{1ex}

Soit $\cal F$ le feuilletage d\'efini par la submersion
\( \pi\colon M \rightarrow B \). Soit $\phi^p$ le faisceau des germes de
p-formes basiques. D'apr\`es \cite{Va}, on sait que
\( H^{q}(M;\phi^{p})  \cong  E_1^{p,q}({\cal F}) \).
On a alors le corollaire suivant qui
g\'en\'eralise le r\'esultat de cohomologie feuillet\'ee de \cite{DH}:

\begin{cor}
Si les fibres $F_b$ sont
connexes et les groupes de cohomologie de De Rham
$H^q(F_b)=0$ en degr\'{e} $q < r$, alors on a
$$
E_1^{p,q}({\cal F}) = \left\{
\begin{array}{ll}
\Omega^p(B) & \mbox{pour $q=0$} \\
0 & \mbox{pour $q < r$}
\end{array}
\right. \;\;\;\;\;\; \Box
$$
\end{cor}

\vspace{4ex}

\begin{flushleft}
Fernando Alcalde Cuesta \\
Departamento de Xeometria e Topoloxia \\
Universidade de Santiago de Compostela  \\
15706 Santiago de Compostela (Espagne) \\
e-mail: alcalde@aimat.cesga.es
\end{flushleft}

\begin{flushleft}
Gilbert Hector \\
Laboratoire de G\'eom\'etrie et Analyse $-$ U.R.A. 746 \\
Universit\'e Claude Bernard-Lyon 1  \\
69622 Villeurbanne (France) \\
e-mail: hector@geometrie.univ-lyon1.fr
\end{flushleft}

\end{document}